\documentclass[a4paper, reqno, 11pt]{amsart}

\usepackage[utf8x]{inputenc}
\usepackage[T1]{fontenc}
\usepackage[ngerman,english]{babel}
\usepackage{amsmath, amsfonts, amsthm, amssymb, amscd}
\usepackage[]{graphicx, psfrag}
\usepackage{bbm,bm}
\usepackage[below]{placeins}
\usepackage{subfigure}
\usepackage{multirow}

\usepackage[dvips]{hyperref}
\hypersetup{
 pdfauthor = {G\"unther~Gr\"un, Fabian~Klingbeil},
 pdftitle = {Two-phase flow with mass density contrast: stable schemes for a
   thermodynamic consistent and frame-indifferent diffuse-interface model}
}

\newcommand{\iOmegaOT}{{\int_0^T \!\!\! \int_\Omega}}
\newcommand{\iOmega}{\int_{\Omega}}
\newcommand{\R}{\mathbb{R}}
\newcommand{\N}{\mathbb{N}}

\newcommand{\abs}[1]{\left|{#1}\right|}
\newcommand{\set}[1]{{\left\{#1\right\}}}
\newcommand{\rkla}[1]{{\left(#1\right)}}

\newcommand{\skla}[1]{{\left\langle#1\right\rangle}}
\newcommand{\ekla}[1]{{\left[#1\right]}}

\newcommand{\vflow}{\mathbf{v}}
\newcommand{\wflow}{\mathbf{w}}
\newcommand{\jflow}{\mathbf{j}}
\newcommand{\Jflow}{\mathbf{J}}
\newcommand{\kflow}{\mathbf{k}}

\newcommand{\Bflow}{\mathbf{B}}
\newcommand{\Dflow}{\mathbf{D}}

\newcommand{\Aflow}{\mathbf{A}}
\newcommand{\Mflow}{\mathbf{M}}
\newcommand{\Nflow}{\mathbf{N}}
\newcommand{\Vflow}{\mathbf{V}}
\newcommand{\Xflow}{\mathbf{X}}

\newcommand{\Th}{\mathcal T_h}
\newcommand{\Thh}{\mathcal T_{\frac{h}{2}}}
\newcommand{\Ph}{\mathcal P_h}
\newcommand{\Jh}{\mathcal I_h}
\newcommand{\Jhh}{\mathcal I_{\frac{h}{2}}}
\newcommand{\Ih}{\mathcal I_h}
\newcommand{\Ihh}{\mathcal I_{\frac{h}{2}}}
\newcommand{\Uh}{U_h}
\newcommand{\Uhh}{U_{\frac{h}{2}}}

\newcommand{\inv}{^{-1}}

\newcommand{\ksemi}{{k+\frac{1}{2}}}
\newcommand{\kimp}{{k+1}}

\newcommand{\rhoDerivative}{\tfrac{\partial\rho}{\partial\varphi}}
\newcommand{\rhoDerivativeKimp}{\tfrac{\delta\rho}{\delta\varphi}}
\renewcommand{\div}{\operatorname{div}}

\newtheorem{theorem}{Theorem}[section]

\newtheorem{corollary}[theorem]{Corollary}

\newtheorem{remark}[theorem]{Remark}

\numberwithin{equation}{section}

\begin{document}
\selectlanguage{english}


\title[Two-phase flow with mass density contrast]{Two-phase flow with mass density contrast: stable schemes for a thermodynamic consistent and frame-indifferent diffuse-interface model}
\author[G. Gr\"un]{G\"unther Gr\"un}
\email{gruen@am.uni-erlangen.de}
\author[F.~Klingbeil]{Fabian Klingbeil}
\address{\selectlanguage{ngerman}G\"unther~Gr\"un, Fabian~Klingbeil: Friedrich-Alexander-Universit\"at Er\-lan\-gen-N\"urn\-berg --- Department Mathematik  --- Cauerstr.~11 --- 91058 Erlangen --- Deutschland}


\date{\today}

\selectlanguage{english}

\begin{abstract}
In this paper, we present a numerical scheme for the diffuse-interface model
in [Abels, Garcke, Gr{\"u}n, M3AS 22(3), 2012] for two-phase flow of immiscible, incompressible fluids. As
that model is in particular consistent with thermodynamics, energy estimates
are expected to carry over to the discrete setting. By a subtle discretization
of the convective coupling with the flux of the phase-field in the momentum
equation, we prove discrete consistency with
thermodynamics. Numerical experiments in two spatial dimensions -- ranging from  Rayleigh-Taylor
instability to a comparison with previous modeling approaches -- indicate the full
practicality of our scheme and enable a first validation of the new modeling
approach in [Abels, Garcke, Gr{\"u}n, M3AS 22(3), 2012].  
\end{abstract}

\maketitle


\section{Introduction}\label{sec:introduction}
 In \cite{AGG}, the following diffuse-interface model for two-phase flow of immiscible, incompressible fluids was introduced.
 \begin{subequations}\label{eq:cont}\begin{gather}
 \rho \partial_t \vflow + \rkla{\rkla{\rho\vflow + \rhoDerivative \jflow} \cdot \nabla} \vflow - \nabla\cdot \rkla{2\eta(\varphi)\Dflow\vflow} + \nabla p = \mu \nabla \varphi + \kflow_\text{grav} \label{eq:cont1} \\
 \partial_t \varphi + \vflow \cdot \nabla \varphi - \nabla\cdot \rkla{M(\varphi) \nabla \mu} = 0 \label{eq:cont2} \\
 \mu = \sigma \rkla{- \Delta \varphi + F'(\varphi)} \label{eq:cont3} \\
 \nabla\cdot \vflow = 0 \quad \text{ in } \Omega \times (0,T) \label{eq:cont4}
\end{gather}\end{subequations}
As boundary conditions, no-slip conditions for $\vflow$ and the vanishing of the normal derivative of $\varphi$ and of $\mu$ on $\partial \Omega \times (0,T)$ are imposed.

Note that system \eqref{eq:cont} constitutes a coupling of a hydrodynamic momentum equation with a Cahn-Hilliard type phase-field equation. $F$ is a double-well potential with minima in $\pm 1$ - representing the pure phases $\varphi\equiv\pm 1$. The parameter $\sigma$ is the surface tension coefficient, which is assumed to be $\sigma=1$ if not stated explicitly. The term $\mu$ stands for the so called chemical potential, and the order parameter $\varphi$ stands for the difference of the volume fractions $u_2 - u_1$ where $u_i(x,t) := \frac{\rho_i(x,t)}{\hat\rho_i}$ with $\hat\rho_i$ the specific (constant) density of fluid $i$ in a unmixed setting. Denoting the individual velocities by $\vflow_i$, $i=1,2$, we write $\vflow:=u_1\vflow_1+u_2\vflow_2$ for the volume averaged velocity. 

The mass density $\rho(\varphi)$ is defined as
\begin{equation*} \rho(\varphi) = \frac{\hat\rho_2 + \hat\rho_1}{2} + \frac{\hat\rho_2 - \hat\rho_1}{2}\varphi \end{equation*}
and $\Dflow\vflow$ denotes the symmetrized gradient. The term $\kflow_\text{grav}$ stands for the density of external volume forces, within this article we only consider gravitational forces. Finally, the flux $\jflow$ is defined by $\jflow := -M(\varphi)\nabla \mu$. 

System \eqref{eq:cont} may be derived from energy considerations and Onsager's variational principle. It is worth mentioning that material frame-indifference, i.e. the invariance of the system with respect to Galilean transformations, for instance, is sensitive to the question to which extent motions relatively to the barycentric motion are neglected in the momentum equation. In the predecessor paper \cite{Abels2010}, relative motions are totally neglected in the momentum equation. As a result, the system studied there turns out not to be frame-indifferent anymore. For more details on the derivation of the model studied here, we refer the reader to \cite{AGG}.
\newline\newline
Conceptually, various approaches exist to model the flow of two immiscible, incompressible fluids. In sharp-interface models, the transition layer separating the two fluids is idealized to be a two-dimensional surface. In level-set methods, volume-of-fluids methods, and diffuse-interface models, additional order parameters are introduced which provide information whether fluid 1 or fluid 2 prevails in a spatial point $x$ at time $t$. 

In contrast to the former approaches, diffuse-interface models assume the transition layer to be of finite size. They are distinguished by the following features. No artificial additional conditions are necessary to model topology changes or to guarantee conservation of individual masses.
Often existence of solutions to the underlying system of partial differential equations can be proven global in time. This is desirable not only from the mathematical point of view. Such a result is available (or expected to be available) in the continuous setting as soon as the model is consistent with thermodynamics.  In our context, this means that the energy at a time $t_2>t_1$ is bounded by the sum of the energy at  time $t_1$ and the work done by external forces during the time interval $[t_1, t_2]$. Especially in the framework of the system under consideration, we require that the sum of the kinetic energy $\frac12 \iOmega \rho \abs{\vflow}^2$ and the interfacial energy $\iOmega \frac12 \abs{\nabla \varphi}^2+F(\varphi)$ is decreasing with respect to time when no external forces are considered. 

Numerically, such a result is the key to prove stability and convergence results for appropriate schemes. In this sense, we introduce the {\it total energy} at time $t$ as the sum of the discrete counterparts of the kinetic and the interfacial energy. We call a numerical scheme {\it stable} or {\it discretely consistent with thermodynamics} if the total energy at time $t_2>t_1$ is bounded by the sum of the total energy at time $t_1$ and the work done by external forces during the time interval $[t_1, t_2]$.

Note that this notion of stability is in the spirit of \emph{stability with respect to a norm}, see for instance \cite[Def. 2.4.3]{Kroner1997}.
\newline\newline
Let us make a few remarks on diffuse-interface models in general.

Historically, the first diffuse-interface model for two-phase flow was the so called model H, introduced by Hohenberg and Halperin \cite{Hohenberg1977}. It assumed the mass densities of the two fluids to be identical. 

Various models were proposed to extend model H also to the case of mass
density contrast (see \cite{Anderson1998} and the references
therein). Lowengrub and Truskinovsky proposed in \cite{Lowengrub1998} for the
first time a diffuse-interface model consistent with thermodynamics. The gross
velocity field is obtained by mass averaging of individual velocities. As a
consequence, it is not divergence free, and the pressure $p$ enters the model
as an essential unknown. However, no energy estimates are available to control
$p$. Moreover, the pressure enters the chemical potential and is hence
strongly coupled to the phase-field equation. This intricate coupling may be
one reason why so far it has not been possible to formulate numerical schemes
for model \cite{Lowengrub1998}. For an existence result, we refer to \cite{Abels2009}. 

Ding et al. \cite{DSS07} suggested  to define the gross velocity field by volume averaging. Prohibiting in addition volume changes due to mixing ("simple mixture assumption"), the gross velocity field is solenoidal. To the best of our knowledge, however, all attempts failed to establish energy inequalities and  to show that the model in \cite{DSS07} is  consistent with thermodynamics.

In \cite{Shen2010}, Shen and Yang propose an extension of the model \cite{DSS07} which allows for energy estimates. Their modeling ansatz is to add a multiple of the term $\rho_t + \div (\rho\vflow)$ in the momentum equation. They justify this idea by the assertion that the continuity equation $\rho_t + \div (\rho\vflow)=0$ were valid and therefore this term were zero. Nevertheless, the phase-field equation $\varphi_t + \div (\varphi\vflow) - \div \Jflow= 0$ is also part of their model, and $\rho$ depends in an affine-linear way on $\varphi$. See also \cite{Salgado2011} for further studies based on the ideas of \cite{Shen2010}. 

A third strategy was pursued by Boyer \cite{Franck2002}, allowing also for solenoidal vector fields, but apparently not for energy estimates. 

In the present  paper, we formulate a fully practical numerical scheme for the model in \cite{AGG} which is discretely consistent with thermodynamics in the sense that energy estimates are satisfied by discrete solutions. 

For the ease of presentation, we present the algorithm in a finite element
context. In section \ref{sec:formulation}, we briefly discuss the model in the
continuous setting and focus on the derivation of the energy estimate. That
way, we motivate our numerical approach, in particular the choice of
projection terms and a new and subtle discretization of the convective term in
the momentum equation. For this scheme, we prove the decay of the discrete energy. This way, we show stability of the scheme.
 
Next, we discuss questions of implementation. It turns out that the scheme proposed in section \ref{sec:formulation}, in particular certain $L^2$-projections, lead to a fully occupied system matrix in the momentum equation. Therefore, we propose a new, equivalent formulation, that avoids the computation of $L^2$-projections. For that scheme, we suggest a splitting ansatz which combines Taylor-Hood-elements for the momentum equation, linear finite elements for phase-field and chemical potential with a finite volume approach to deal with the convective terms. More precisely, we discretize transport terms on the dual mesh associated with the finite element grid and use second-order methods for scalar conservation laws, i.e. upwind schemes with min-mod limiters (see \cite{Kroner1997}). 


In the last section, we present characteristic two-dimensional numerical experiments to underline the good performance of our scheme. Topics include rising bubbles and the simulation of Rayleigh-Taylor instabilities. We allow for Atwood numbers $A := \frac{\hat\rho_1 - \hat\rho_2}{\hat\rho_1 + \hat\rho_2}$ up to 0.99 and we investigate the effect caused by significant differences in the viscosities. 

Finally, we present numerical experiments for annulus-shaped droplets in oscillating force-fields which show considerable differences between the results obtained based on the model of Ding et al. \cite{DSS07} and those based on \cite{AGG}.  

\section{A stable  scheme}\label{sec:formulation}

\subsection{Basic assumptions}\label{subsec:preleminaries}
We consider the two-phase problem on a bounded, convex polygonal (or polyhedral, respectively) domain $\Omega\subset\R^d$ in spatial dimensions $d\in\{2,3\}$.  By $\skla{ \cdot, \cdot }$, we denote the Euclidean scalar product on $\R^d$. We assume $\Th$ to be a regular and admissible triangulation of $\Omega$ with simplicial elements in the sense of \cite{Ciarlet2002}.
By $\Uh$, we denote the space of continuous, piecewise linear finite element functions on $\Th$. We write $\Thh$ for the regular and admissible triangulation of $\Omega$ which is obtained by introducing the mid-points of edges in $\Th$ as additional nodes. $\Uhh$ denotes the space of continuous, piecewise linear finite element functions on $\Thh$. The expression $\Jh$ (or $\Jhh$, respectively) stands for the nodal interpolation operator from $C^0(\Omega)$ to $\Uh$ (or $\Uhh$, respectively) defined by $\Jh u:=\sum_{j=1}^{\dim \Uh }u(x_j)\psi_j$, where the functions $\psi_j$ form a dual basis to the nodes $x_j$, i.e. $\psi_i(x_j) = \delta_{ij}$, $i,j=1, \dots, \dim \Uh$. The definition of $\Jhh$ is completely analogous. By $\Ph$, we denote the orthogonal $L^2-$projection onto $\Uh$. Moreover, we introduce the finite element spaces
\begin{align*}
  \Uh^0 & := \set{\psi \in C^0(\bar\Omega): \psi|_K \in P_1(K), K \in \Th ,\int_\Omega \psi \equiv 0}, \\ 
 \Xflow_h & := \set{\wflow \in (C^0_0(\bar\Omega))^d : \rkla{\wflow}_j|_K \in  P_2(K), K \in \Th, j=1,\dots,d }, \quad d=2,3,\\
 \Vflow_h & := \set{\wflow \in \Xflow_h: \iOmega  \psi \div \wflow =0 \, \forall \, \psi \in \Uh^0}. 
\end{align*}
and write $\set{\wflow_1,\dots, \wflow_{\dim \Xflow_h}}$ to denote the dual basis on $\Xflow_h$.

Note that the degrees of freedom used for each component of $\Xflow_h$ on
an arbitrary element $K\in\Th$ are in agreement with the degrees of freedom of
elements of $\Uhh$ restricted to $K$. The reader may change the space
$\Xflow_h$ in a straightforward way to adapt to a flow with a Navier-slip
boundary condition $\vflow_\tau=\beta\partial_n\vflow_\tau$.
Concerning the discretization with respect to time, we consider a time interval $[0,T]$ and subdivide it for $N\in\N$ into $N$ subintervals $I_k:=[t_{k-1}, t_k]$, $k=1,\dots,N$, of length $\tau:=\frac T N$. Non-constant time-increments are also possible -- for the ease of presentation we confine ourselves to the uniform case. Note that for a function $f:[0,T] \to X$, we abbreviate $f^k:=f(t_k)$, $k=0,\dots,N$.

\subsection{The model in the continuous setting}\label{subsec:continousmodel}

In this section, we sketch the proof of an energy estimate for \eqref{eq:cont} in the continuous setting. This way, we intend to motivate the particular discretization to be proposed in the next subsection.  For the ease of presentation, we assume external forces here to vanish identically.

Formally, we have the following energy identity.
\begin{multline}\label{eq:cont_energyest}
 \tfrac12 \iOmega \rho(T) \abs{\vflow}^2(T) + \tfrac12 \iOmega \abs{\nabla
   \varphi}^2(T) + \iOmega F(\varphi(\cdot, T))  
  + \iOmegaOT M(\varphi)\abs{\nabla \mu}^2 \\ + \iOmegaOT
  2\eta(\varphi)\abs{\Dflow\vflow}^2 = \tfrac12 \iOmega \rho^0
  \abs{\vflow^0}^2 + \tfrac12 \iOmega \abs{\nabla \varphi^0}^2 
  + \iOmega F(\varphi^0)
\end{multline}
where $\varphi^0(\cdot) := \varphi(\cdot, 0)$.

Indeed,  multiplying \eqref{eq:cont1}, \eqref{eq:cont2} and \eqref{eq:cont3} by $\vflow$, $\mu + \frac12\rhoDerivative\abs{\vflow}^2$ and $\partial_t \varphi$, respectively, we obtain the  identities
\begin{subequations}\label{estimates}
\begin{multline} \label{estimates1}
\tfrac12 \iOmega \rho\partial_t|\vflow|^2- \tfrac12 \iOmega \rho'\skla{\vflow,\nabla\varphi}|\vflow|^2 + \tfrac12 \iOmega \rho'\skla{\jflow,\nabla|\vflow|^2}  \\
 + \iOmega \eta(\varphi)|\Dflow\vflow|^2= \iOmega \mu \skla{\vflow,\nabla\varphi}
\end{multline}
\begin{equation} \label{estimates2}
\tfrac12 \partial_t \iOmega|\nabla\varphi|^2+\partial_t \iOmega F(\varphi) + \iOmega M(\varphi)|\nabla\mu|^2= -\iOmega \mu \skla{\vflow,\nabla\varphi} 
\end{equation}
\begin{equation} \label{estimates3}
\tfrac12\iOmega\partial_t\rho|\vflow|^2 + \tfrac12 \iOmega \rho'\skla{\vflow,\nabla\varphi}|\vflow|^2 - \tfrac12 \iOmega \rho'\skla{\jflow,\nabla|\vflow|^2}=0. 
\end{equation}
\end{subequations}
Summing up and integrating with respect to time, the energy identity \eqref{eq:cont_energyest} follows.

\subsection[A mixed finite element scheme]{A mixed finite element scheme --
  discretely consistent with thermodynamics}\label{subsec:fe-methods} 

In this section, we propose a finite element scheme which allows for a discrete counterpart of the energy estimate \eqref{eq:cont_energyest} and which is therefore consistent with thermodynamics in a discrete sense. 

For the velocity field and the pressure, we use Taylor-Hood finite elements, the phase-field (and consequently the mass density) will be approximated by linear finite elements from $\Uh$. The same holds for the chemical potential $\mu$.

We introduce the notation 
\begin{equation}\label{eq:diffquotient}
\rhoDerivativeKimp:=\begin{cases}
\tfrac{\rho(\varphi^\kimp)-\rho(\varphi^k)}{\varphi^\kimp-\varphi^k}
& \text{ if } \varphi^\kimp\ne\varphi^k\\
\tfrac{\partial\rho(\varphi^\kimp)}{\partial\varphi} & \text{ otherwise.}
\end{cases}
\end{equation}
The motivation is as follows. For the energy estimate, we need $\rho(\varphi)$ to be bounded away from zero. Therefore, we replace $\rho(\varphi)$ by an appropriate (nonlinear) regularization $\bar\rho(\varphi)$, establish the energy estimate and choose the double-well potential $F(\cdot)$ depending on the grid size $h$ in such a way that $\varphi$ stays in $(-1-\epsilon, 1+\epsilon)$ and that way in the linear regime of $\rho$. 

Let us now present our strategy to guarantee discrete consistency with thermodynamics. Observe that the derivation of \eqref{eq:cont_energyest} requires testing the phase-field equation \eqref{eq:cont2} by a multiple of $\abs{\vflow}^2$ which cannot be expected to be an admissible test function in the discrete setting. 

Furthermore, a closer look at \eqref{estimates1} and \eqref{estimates3} reveals that the product of $\vflow$ with the convective terms in \eqref{eq:cont1} cancels out against the product of $\frac12 \rhoDerivative \abs{\vflow}^2$ with the second and third term in \eqref{estimates3}. 

Hence, it seems to be crucial to discretize the convective term in \eqref{eq:cont1} in such a way that the cancellation carries over to the discrete setting. Similarly, the first term in \eqref{eq:cont1} requires a special time averaging of the mass density to guarantee a control of the kinetic energy in the discrete setting. 

Here is our approach to discretize the convective terms in \eqref{eq:cont1}. Starting from the identity 
\begin{multline*}
 \iOmega \rkla{(\jflow \cdot \nabla)\vflow} \cdot \wflow = \tfrac12 \iOmega \skla{\jflow , \nabla \skla{\vflow, \wflow}} \\
 + \tfrac12\iOmega \skla{\jflow , \rkla{(\nabla\vflow)^T \wflow}} - \tfrac12 \iOmega \skla{\jflow , \rkla{(\nabla\wflow)^T \vflow}}
\end{multline*}
we discretize $\iOmega \rhoDerivative \rkla{(\jflow \cdot \nabla)\vflow} \cdot
\wflow$ by 
\begin{multline}\label{eq:disc_rhoJ}
 \iOmega \rhoDerivative \rkla{(\jflow \cdot \nabla)\vflow} \cdot \wflow
  \approx \tfrac12 \iOmega  \skla{\jflow^\kimp, \nabla \Ph\rhoDerivativeKimp \Jhh\skla{\vflow^k, \wflow}} \\
   + \tfrac12 \iOmega \rhoDerivativeKimp \skla{\jflow^\kimp, (\nabla\vflow^\kimp)^T \wflow}
  - \tfrac12 \iOmega \rhoDerivativeKimp \skla{\jflow^\kimp, (\nabla\wflow)^T \vflow^\kimp}.
\end{multline}
Similarly, for solenoidal vector fields $\vflow$ with non-penetrating boundary
data ($\vflow\cdot\mathbf n=0$), the identity 
\begin{multline*}
 \iOmega \rho \rkla{ \rkla{\vflow\cdot\nabla} \vflow}\wflow
  = -\tfrac12 \iOmega \rhoDerivative \skla{\vflow, \nabla\varphi} \skla{\vflow, \wflow} \\
   - \tfrac12 \iOmega \rho \skla{\vflow, \rkla{\nabla\wflow}^T \vflow} + \tfrac12 \iOmega \rho \skla{\vflow, \rkla{\nabla\vflow}^T \wflow}
\end{multline*}
suggests to discretize $\iOmega \rho \rkla{\rkla{\vflow \cdot \nabla}\vflow}\cdot \wflow$ by
\begin{multline}\label{eq:disc_rhoVelo}
 \iOmega \rho \rkla{\rkla{\vflow \cdot \nabla}\vflow}\cdot \wflow
  \approx -\tfrac12 \iOmega \skla{\vflow^\kimp, \nabla\varphi^\kimp} \Ph \rhoDerivativeKimp \Jhh{\skla{\vflow^k, \wflow}} \\
   - \tfrac12 \iOmega \rho^k \skla{\vflow^k, (\nabla\wflow)^T \vflow^\kimp} + \tfrac12 \iOmega \rho^k \skla{\vflow^k, (\nabla\vflow^\kimp)^T \wflow}.
\end{multline}
The question remains how to discretize the mass density factor in the first term of the momentum equation (\ref{eq:cont1}). It turns out that time averaging is the right approach -- see the subsequent Theorem~\ref{lemma:energy}. Summing up, the new scheme reads as follows:

\paragraph{\bf Scheme A}
For given $\varphi^0 \in \Uh$, $\vflow^0 \in \Vflow_h$, $T>0$ and $k=0, \dots, N-1$, find functions $\rkla{\varphi^\kimp, \mu^\kimp, \vflow^\kimp} \in \Uh \times \Uh \times \Vflow_h$ such that with $\tau := \frac{T}{N}$ we have
\begin{subequations}\label{eq:scheme}
\begin{multline}\label{eq:scheme_a}
 \iOmega \frac{\rho^k + \rho^\kimp}{2} \Jhh \skla{\partial^-_\tau \vflow^\kimp, \wflow}
  - \tfrac12 \iOmega \skla{\vflow^\kimp, \nabla\varphi^\kimp} \Ph \rhoDerivativeKimp \Jhh{\skla{\vflow^k, \wflow}} \\
  - \tfrac12 \iOmega \rho^k \skla{\vflow^k, (\nabla\wflow)^T \vflow^\kimp}
  + \tfrac12 \iOmega \rho^k \skla{\vflow^k, (\nabla\vflow^\kimp)^T \wflow} \\
  + \tfrac12 \iOmega \skla{\jflow^\kimp, \nabla \Ph\rhoDerivativeKimp\Jhh\skla{\vflow^k, \wflow}}
  + \tfrac12 \iOmega \rhoDerivativeKimp \skla{\jflow^\kimp, (\nabla\vflow^\kimp)^T \wflow} \\
  - \tfrac12 \iOmega \rhoDerivativeKimp \skla{\jflow^\kimp, (\nabla\wflow)^T \vflow^\kimp}
  + \iOmega 2\eta^k \Dflow\vflow^\kimp : \Dflow\wflow  \\
  = \iOmega \mu^\kimp \skla{\nabla\varphi^\kimp, \wflow} +\iOmega \skla{\kflow_\text{grav}(t_k) , \wflow} \quad \text{ for all }\wflow \in \Vflow_h.
  \end{multline}
\begin{multline}\label{eq:scheme_b}
 \iOmega \partial^-_\tau \varphi^\kimp \psi + \iOmega\skla{\vflow^\kimp, \nabla\varphi^\kimp}\psi \\
  + \iOmega M(\varphi^k) \skla{\nabla\mu^\kimp, \nabla\psi} = 0 \quad \text{ for all }\psi \in \Uh
\end{multline}
\begin{multline}\label{eq:scheme_c}
 \iOmega \mu^\kimp \psi = \iOmega \skla{\nabla \varphi^\kimp, \nabla\psi} \\
  + \iOmega \Jh \rkla{\rkla{F'_+(\varphi^\kimp) + F'_-(\varphi^k)} \psi} \quad \text{ for all }\psi \in \Uh
\end{multline}
\end{subequations}
Here, we abbreviated $\partial^-_\tau \Psi^\kimp := \frac{1}{\tau} \rkla{\Psi^\kimp - \Psi^k}$ and $\jflow^\kimp := -M(\varphi^k) \nabla \mu^\kimp$. With $F_+$ and $F_-$, we denote the convex (or concave, respectively) part of $F=F_+ + F_-$. In addition, $\rho^k := \rho(\varphi^k)$ and $\eta^k := \eta(\varphi^k)$. Recall that $\rhoDerivativeKimp = \frac{\rho^\kimp-\rho^k}{\varphi^\kimp - \varphi^k}$. Let us mention already here, that in the course of practical computations we never observed any necessity to apply a nonlinear approximation of $\rho(\cdot)$.   

Note that the nodal interpolation operator $\Ihh$ may be used to simplify the computation of $\Ph \rhoDerivativeKimp \skla{\vflow^k, \wflow}$, see the discussion in Subsection \ref{subsec:reformulation}. In the same subsection, we also present a different, but mathematically equivalent version of \eqref{eq:scheme} which does not involve projection operators. For that scheme, the nodal interpolation operator might be omitted.

Let us now formulate the announced result on stability of the scheme -- or synonymously -- on discrete consistency with thermodynamics.

\begin{theorem}\label{lemma:energy}
Assume that the triple $\rkla{\varphi^\kimp, \mu^\kimp, \vflow^\kimp}$ solves \eqref{eq:scheme} for given $\rkla{\varphi^k, \mu^k, \vflow^k}$. Then
\begin{multline*}
\frac{1}{2\tau} \ekla{ \iOmega \rho^\kimp \Jhh\abs{\vflow^\kimp}^2 - \iOmega \rho^k \Jhh\abs{\vflow^k}^2 + \iOmega \rho^k \Jhh\abs{\vflow^\kimp - \vflow^k}^2 } \\
 + \frac{1}{2\tau} \ekla{ \iOmega \abs{\nabla\varphi^\kimp}^2 - \iOmega \abs{\nabla\varphi^k}^2 + \iOmega \abs{\nabla\varphi^\kimp - \nabla\varphi^k}^2 } \\
  + \frac{1}{\tau} \iOmega \Jh\rkla{F(\varphi^\kimp) - F(\varphi^k)} + \iOmega M(\varphi^k) \abs{\nabla\mu^\kimp}^2 + \iOmega 2\eta^k \abs{\Dflow\vflow^\kimp}^2 \\
\leq \iOmega \skla{\kflow_\text{grav}(t_k), \vflow^\kimp}
\end{multline*}
\end{theorem}

One easily establishes the following global version of the stability estimate which precisely states the boundedness of the total energy in terms of the sum of initial energy and work done by external forces. 

\begin{corollary}
Assume the time interval $[0,T]$ to be subdivided in subintervals $[t_{k-1}, t_k]$ corresponding to time-steps $t_k=k\tau$, $k=0,\dots ,N$, with a uniform time-increment $\tau:=\frac T N$ (for the ease of presentation). With the notation
\begin{equation*}
\mathcal E(\rho,\varphi,\vflow)(t) := \rkla{\tfrac12\iOmega\rho\Jhh|\vflow|^2 + \tfrac12 \iOmega|\nabla\varphi|^2 + \iOmega \Ih \rkla{F(\varphi)} }(t),
\end{equation*}
we have
\begin{multline*}
 \mathcal E(\rho,\varphi,\vflow)(t_k)\\
+ \tau\sum_{m=l}^{k-1} \iOmega \rho(t_m) \Jhh \abs{\vflow(t_{m+1})-\vflow(t_m)}^2 +
\frac\tau 2 \sum_{m=l}^{k-1} \iOmega \abs{\nabla\varphi(t_{m+1})-\nabla\varphi(t_m)}^2 \\
+ \tau \sum_{m=l}^{k-1} \iOmega M(\varphi(t_m)) \abs{\nabla\mu(t_{m+1})}^2 
+ \tau \sum_{m=l}^{k-1} \iOmega 2 \eta(\varphi(t_m)) \abs{\Dflow\vflow(t_{m+1})}^2 \\
 \leq \mathcal E(\rho,\varphi,\vflow)(t_l) + \tau\sum_{m=l}^{k-1} \iOmega \skla{\kflow_\text{grav}(t_{m}), \vflow(t_{m+1})}
\end{multline*}
for any $0\leq l<k\leq N$ in $\N$.
\end{corollary} 

\begin{remark}
 \begin{enumerate}
  \item
   Extending the double-well potential for each $h>0$ in an appropriate way into the complement of $[-1, 1]$, it is possible to bound $\varphi$ to attain values in the interval $(-1-\varepsilon, 1+\varepsilon)$ for given $\varepsilon$. 
  \item
By using mathematical fixed-points results, it is possible to prove existence of discrete solutions.
 Details will be presented elsewhere.
   
 \end{enumerate}
\end{remark}

\begin{proof}
Choosing $\psi := \partial^-_\tau \varphi^\kimp$ in \eqref{eq:scheme_c} and $\psi := \mu^\kimp$ in \eqref{eq:scheme_b} yields -- using the convexity of $F_+$ and $-F_-$ -- after summation
\begin{multline}\label{eq:energy_c}
 \frac{1}{2\tau} \ekla{ \iOmega \abs{\nabla\varphi^\kimp}^2 - \iOmega \abs{\nabla\varphi^k}^2 + \iOmega \abs{\nabla\varphi^\kimp - \nabla\varphi^k}^2 } + \iOmega M(\varphi^k)\abs{\nabla\mu^\kimp}^2 \\
  + \frac{1}{\tau} \iOmega \Jh\rkla{F(\varphi^\kimp) - F(\varphi^k)} +  \iOmega \skla{\vflow^\kimp, \nabla\varphi^\kimp}\mu^\kimp  \leq 0
\end{multline}
 Testing \eqref{eq:scheme_b} by $\psi := \tfrac12 \Ph \rhoDerivativeKimp \Jhh\skla{\vflow^k, \vflow^\kimp}$ and using the $L^2$-projection property of $\Ph$ yields
 \begin{multline}\label{eq:energy_a}
  \tfrac12 \iOmega \partial^-_\tau \rho^\kimp \Jhh\skla{\vflow^k, \vflow^\kimp} + \tfrac12\iOmega \skla{\vflow^\kimp, \nabla\varphi^\kimp} \Ph \rhoDerivativeKimp \Jhh\skla{\vflow^k, \vflow^\kimp} \\
   - \tfrac12 \iOmega  \skla{\jflow^\kimp, \nabla \Ph \rhoDerivativeKimp\Jhh\skla{\vflow^k, \vflow^\kimp}} = 0
\end{multline}
 Testing \eqref{eq:scheme_a} by $\wflow := \vflow^\kimp$ and using the identity
 \begin{multline*}
  \rkla{\rho^k + \rho^\kimp} \Jhh\skla{\vflow^\kimp - \vflow^k, \vflow^\kimp} 
   = \rho^\kimp \Jhh\abs{\vflow^\kimp}^2 - \rho^k \Jhh\abs{\vflow^k}^2 \\
   + \rho^k \Jhh\abs{\vflow^\kimp - \vflow^k}^2 - (\rho^\kimp - \rho^k) \Jhh\skla{\vflow^k, \vflow^\kimp}
 \end{multline*}
 entails
 \begin{multline}\label{eq:energyb_b}
  \frac{1}{2\tau} \ekla{ \iOmega\rho^\kimp \Jhh\abs{\vflow^\kimp}^2 - \iOmega\rho^k\Jhh\abs{\vflow^k}^2 + \iOmega\rho^k \Jhh\abs{\vflow^\kimp - \vflow^k}^2 } \\
  - \tfrac12 \iOmega \partial^-_\tau \rho^\kimp \Jhh\skla{\vflow^k, \vflow^\kimp} - \tfrac12 \iOmega \skla{\vflow^\kimp, \nabla\varphi^\kimp} \Ph \rhoDerivativeKimp \Jhh\skla{\vflow^k, \vflow^\kimp} \\
  + \tfrac12 \iOmega  \skla{\jflow^\kimp, \nabla \Ph\rhoDerivativeKimp\Jhh\skla{\vflow^k, \vflow^\kimp}} + \iOmega 2\eta^k \abs{\Dflow\vflow^\kimp}^2  \\
   = \iOmega \mu^\kimp \skla{\vflow^k, \nabla\varphi^\kimp}+\iOmega \skla{ \kflow_\text{grav}(t_k), \vflow^\kimp}
 \end{multline}
 Adding \eqref{eq:energy_c}, \eqref{eq:energy_a} and \eqref{eq:energyb_b} gives the result.
\end{proof}

\section{Implementation}\label{sec:implementation}

\subsection{Reformulation without projection terms}\label{subsec:reformulation}

In this section, we discuss practical aspects of \eqref{eq:scheme}. A first ansatz to compute discrete solutions would be to split the problem into solving the Cahn-Hilliard part and the momentum equation part separately. A closer look at equation \eqref{eq:scheme_a} reveals that the second term on the left-hand side requires additional approximation, even if $\varphi^\kimp$ is assumed to be known already. 
Indeed, expanding $\vflow^\kimp=\sum_{i=1}^{\dim \Xflow_h}V_i^\kimp\wflow_i$ in terms of a basis of $\Xflow_h$, we calculate
\begin{equation*}
 \iOmega \skla{\vflow^\kimp, \nabla\varphi^\kimp} \Ph\rhoDerivativeKimp\Jhh\skla{\vflow^k, \wflow_i} = \rkla{\mathbf{S}^k \cdot \begin{pmatrix}
 		V^\kimp_1 \\
 		\vdots \\
 		V^\kimp_{\dim\Xflow_h} \\
 \end{pmatrix}}_i,
\end{equation*}
where the matrix $\mathbf{S}^k \in \R^{\rkla{\dim\Xflow_h}^2}$ has components
\begin{equation*}
 \rkla{\mathbf{S}^k}_{i,j} := \iOmega \skla{\wflow_j, \nabla\varphi^\kimp} \Ph\rhoDerivativeKimp\Jhh\skla{\vflow^k, \wflow_i}, \quad i,j=1, \dots, \dim \Xflow_h.
\end{equation*}
Observe that $ \Ph \rhoDerivativeKimp \Jhh\skla{\vflow^k, \wflow_i}$ can be written as
\begin{equation*}
\Ph \rhoDerivativeKimp \Jhh\skla{\vflow^k, \wflow_i} = \sum_{j=1}^{\dim \Uh} \hat p^k_{i,j} \psi_j(x)
\end{equation*}
with
\begin{equation*} 
\hat p_{i}^k := \mathbf{M}_{\Uh}\inv \begin{pmatrix}
		\iOmega \rhoDerivativeKimp \Jhh \skla{\vflow^k, \wflow_i} \psi_1 \\
		\vdots \\
		\iOmega \rhoDerivativeKimp \Jhh \skla{\vflow^k, \wflow_i} \psi_{\dim \Uh}  \\
\end{pmatrix}.
\end{equation*}
Here, $\mathbf{M}_{\Uh}$ denotes the mass matrix $\rkla{\mathbf{M}_{\Uh}}_{ij} = \iOmega \psi_i \psi_j$, the inverse of which
is fully occupied. As a consequence, $\mathbf{S}_{ij}^k \neq 0$ for all $i,j=1,\dots,\dim\Xflow_h$ in
general. Hence, the system matrix may become fully occupied in general,
which seems not to be tolerable with respect to numerical costs.

Of course, it would be possible to take advantage of the fact that $\Jhh \skla{\vflow^k,\wflow}$ is contained in $\Uhh$ and that $\rhoDerivativeKimp$ can be assumed to be constant in practical computations\footnote{All our numerical experiments suggest that in practice it is not necessary to replace $\rho(\varphi)$ by a nonlinear approximation bounded away from zero since $\varphi$ always stayed sufficiently close to $[-1,1]$.}. This would simplify the computation of the projection.

Expanding a basis of $\Uh$ to a basis of $\Uhh$, projections are to be computed only for the basis functions which are not contained in $\Uh$. For these basis functions, an approximate result may be obtained by replacing $\vflow^\kimp$ by $\vflow^k$. This way, the corresponding terms do not enter the system matrix anymore (see \cite{Grun2011}).

In the present manuscript, however, we pursue a different approach which allows to avoid the computation of orthogonal projections at all. We proceed as follows.

Testing equation \eqref{eq:scheme_b} by $\psi=\frac12 \Ph\rhoDerivative \Jhh \skla{\vflow^k,\wflow}$, we obtain the identity
\begin{multline*}
\tfrac12\int_\Omega \partial^-_\tau \rho^{k+1}\Jhh \skla{\vflow^k,\wflow} +
\tfrac12\int_\Omega \skla{\vflow^{k+1},\nabla\varphi^{k+1}} \Ph\rhoDerivative \Jhh \skla{\vflow^k,\wflow} \\
 - \tfrac12 \int_\Omega \skla{\jflow^{k+1},\nabla \Ph\rhoDerivative \Jhh \skla{\vflow^k,\wflow}} =0.
\end{multline*}
Inserting this term into \eqref{eq:scheme_a}, one has
 \begin{align}
\label{eq:scheme_new}
\begin{split}
 \iOmega \frac{\rho^k + \rho^\kimp}{2} \Jhh \skla{\partial^-_\tau \vflow^\kimp, \wflow}
  + \tfrac12 \iOmega \partial_\tau^-\rho^{k+1} \Jhh \skla{\vflow^k,\wflow} \\
  - \tfrac12 \iOmega \rho^k \skla{\vflow^k, (\nabla\wflow)^T \vflow^\kimp}
  + \tfrac12 \iOmega \rho^k \skla{\vflow^k, (\nabla\vflow^\kimp)^T \wflow} \\
   + \tfrac12 \iOmega \rhoDerivativeKimp \skla{\jflow^\kimp, (\nabla\vflow^\kimp)^T \wflow} 
  - \tfrac12 \iOmega \rhoDerivativeKimp \skla{\jflow^\kimp, (\nabla\wflow)^T \vflow^\kimp}\\
  + \iOmega 2\eta^k \Dflow\vflow^\kimp : \Dflow\wflow  \\
  = \iOmega \mu^\kimp \skla{\nabla\varphi^\kimp, \wflow} + \iOmega \skla{\kflow_\text{grav}(t_k), \wflow} \quad \text{ for all }\wflow \in \Vflow_h.
\end{split}
\end{align}
Summing up and reformulating the momentum equation equivalently with test functions which are not discretely divergence free, we obtain the following scheme.

\paragraph{\bf Scheme B} For given $\varphi^0 \in \Uh$, $\vflow^0 \in \Vflow_h$, $T>0$ and $k=0, \dots, N-1$, find functions $\rkla{\varphi^\kimp, \mu^\kimp, \vflow^\kimp, p^\kimp} \in \Uh \times \Uh  \times \Xflow_h \times \Uh^0$ such that with $\tau := \frac{T}{N}$ we have
\begin{subequations}\label{eq:schemeB}
\begin{multline}\label{eq:schemeB_a}
 \iOmega \frac{\rho^k + \rho^\kimp}{2} \Jhh \skla{\partial^-_\tau \vflow^\kimp, \wflow}
  + \tfrac12 \iOmega \partial_\tau^-\rho^{k+1} \Jhh \skla{\vflow^k,\wflow} \\
  - \tfrac12 \iOmega \rho^k \skla{\vflow^k, (\nabla\wflow)^T \vflow^\kimp}
  + \tfrac12 \iOmega \rho^k \skla{\vflow^k, (\nabla\vflow^\kimp)^T \wflow} \\
   + \tfrac12 \iOmega \rhoDerivativeKimp \skla{\jflow^\kimp, (\nabla\vflow^\kimp)^T \wflow} 
  - \tfrac12 \iOmega \rhoDerivativeKimp \skla{\jflow^\kimp, (\nabla\wflow)^T \vflow^\kimp} \\
  + \iOmega 2\eta^k \Dflow\vflow^\kimp : \Dflow\wflow - \iOmega p^\kimp \div \wflow  \\
  = \iOmega \mu^\kimp \skla{ \nabla\varphi^\kimp, \wflow} +\iOmega \skla{\kflow_\text{grav}(t_k), \wflow} \quad \text{ for all } \wflow \in \Xflow_h.
  \end{multline}
  \begin{equation}\label{eq:schemeB_aa}
 \iOmega \psi \div \vflow^\kimp = 0 \quad \text{ for all } \psi \in \Uh^0
\end{equation}
\begin{multline}\label{eq:schemeB_b}
 \iOmega \partial^-_\tau \varphi^\kimp \psi + \iOmega\skla{\vflow^\kimp, \nabla\varphi^\kimp}\psi \\
  + \iOmega M(\varphi^k) \skla{\nabla\mu^\kimp, \nabla\psi} = 0 \quad \text{ for all }\psi \in \Uh
\end{multline}
\begin{multline}\label{eq:schemeB_c}
 \iOmega \mu^\kimp \psi = \iOmega \skla{\nabla \varphi^\kimp, \nabla\psi} \\
  + \iOmega \Jh \rkla{\rkla{F'_+(\varphi^\kimp) + F'_-(\varphi^k)}\psi} \quad \text{ for all }\psi \in \Uh
\end{multline}
\end{subequations}
Note that in equation \eqref{eq:schemeB_a} the nodal interpolation operator $\Ihh$ might be omitted -- to the prize of a slightly higher numerical cost and for a negligible gain in precision.

\subsection{A splitting method}\label{eq:subsec:splitting_method}
In this section we formulate a splitting scheme for \eqref{eq:schemeB} which takes the strong coupling between the momentum equation \eqref{eq:schemeB_a} and the convective Cahn-Hilliard equation \eqref{eq:schemeB_b} and \eqref{eq:schemeB_c} into account and which distinguishes in particular between diffusive and convective parts in the Cahn-Hilliard equation.
We will first describe separately how new iterates for the Cahn-Hilliard equation (Problem CH) and for the momentum equation (Problem ME) are computed. As a third step, we discuss the iterative aspects of our splitting scheme.
As before we write $f^k:=f(t_k)$ for the value of a function $f$ in the $k$-th time-step $t_k$.

\paragraph{{\bf Problem (CH)}}\label{par:problem_ch}\addcontentsline{toc}{subsubsection}{Problem (CH)}
Let us begin with the computation of a new Cahn-Hilliard iterate. We treat the convective part by a finite volume approach and the diffusive part by a finite element approach. We are given a simplicial finite element triangulation $\Th$, adaptively refined by bisection. By $\bar\Th$, we denote the dual grid to $\Th$. In particular, the cell $T_i \in \Th$ dual to a grid point $x_i$ is defined as
\begin{equation*} 
T_i := \set{x \in \Omega: \abs{x-x_i}\leq \abs{x-x_j} \text{ for all } j \neq i, j \in \set{1,..,\dim \Uh}}
\end{equation*} 

 \begin{figure}[ht]
  \centering
  \subfigure{ \includegraphics[width=4cm]{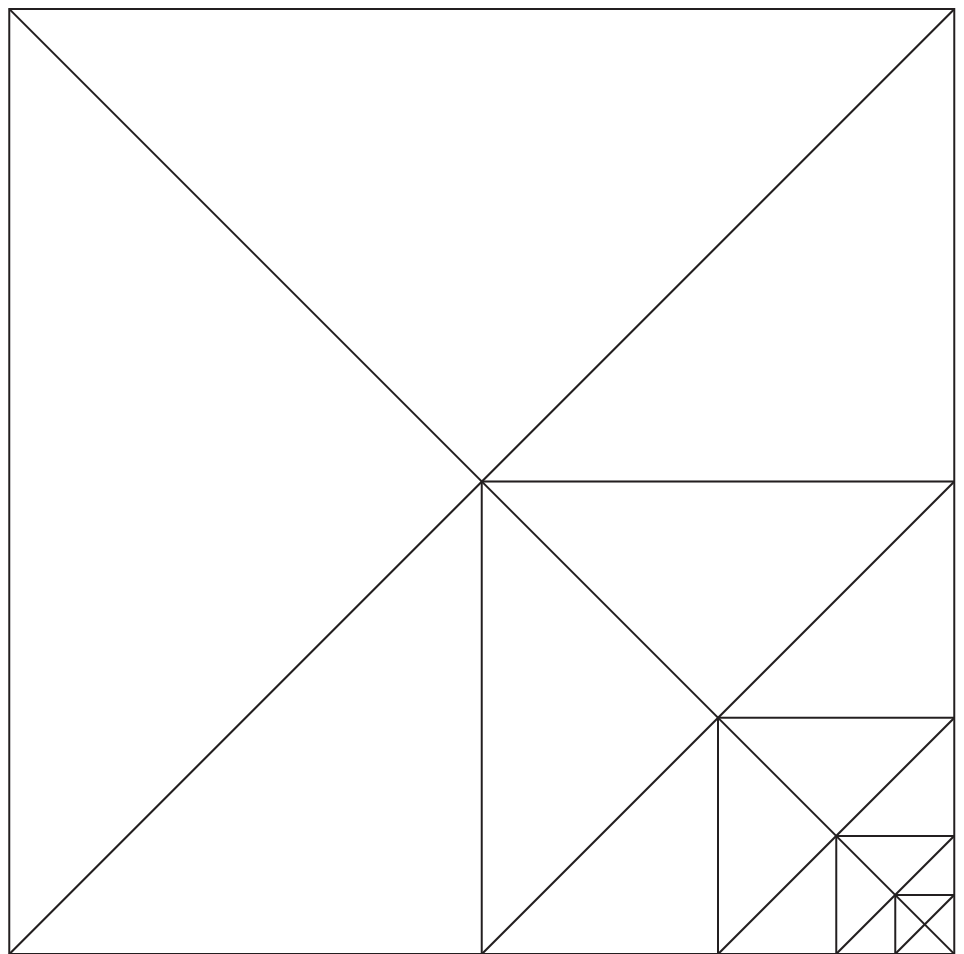} \label{fig:grid_primal} }
  \subfigure{ \includegraphics[width=4cm]{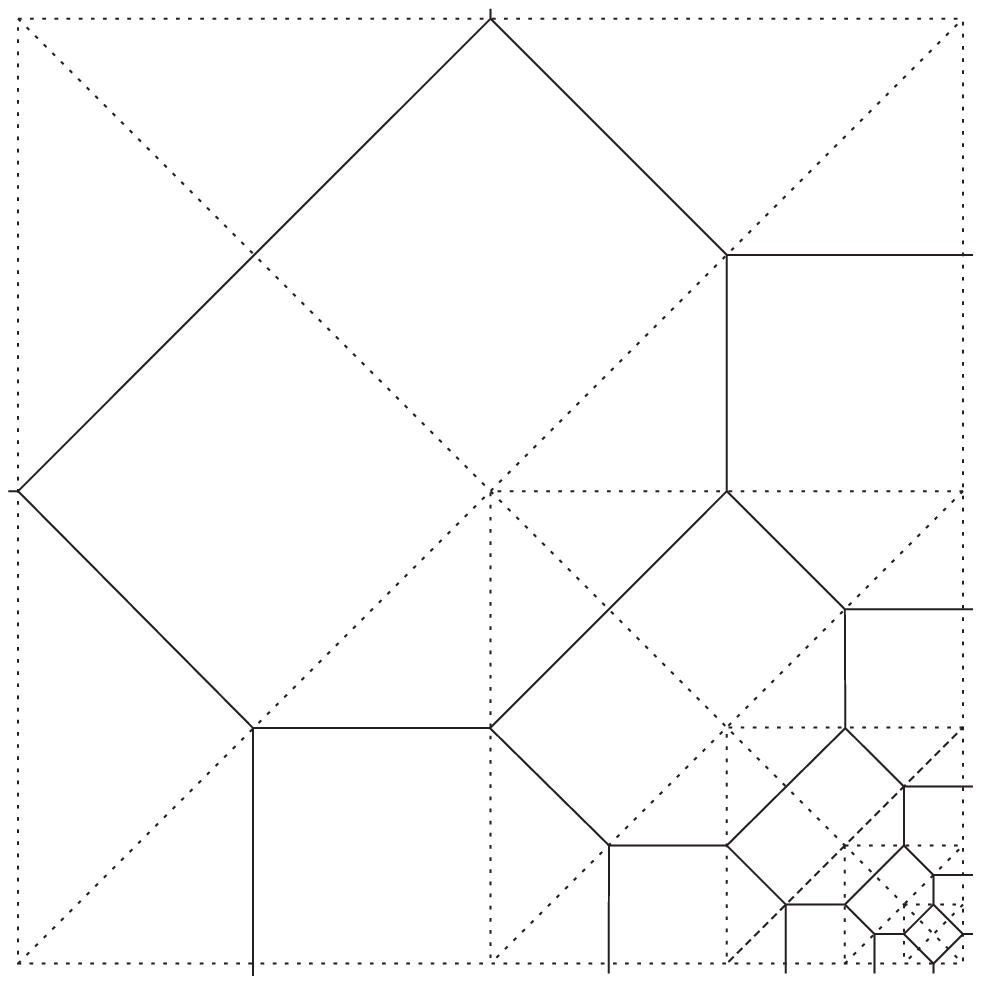} \label{fig:grid_dual} }
  \caption{Primal grid $\Th$ and dual grid $\bar\Th$}
  \label{fig:grid}
 \end{figure}
 
The outer normal $\nu_{ij}$ of the cell $T_i$, pointing towards a neighboring cell $T_j$, is defined by $\nu_{ij} := \frac{x_j -x_i}{\abs{x_j - x_i}}$.
 
For a given finite element function $\varphi \in \Uh$, we denote the corresponding finite volume function by $\bar\varphi$. More precisely, $\bar\varphi$ is elementwise constant in $\bar\Th$, and for all nodal points $x_i, i=1,...,\dim \Uh$, we have $\varphi(x_i) = \bar\varphi(x_i) =: \varphi_i$.

Assume that data $\vflow$ and $\varphi^k$ are given for the velocity field and phase-field, respectively.

Treating the convective and the diffusive part separately, we first use a second order finite volume scheme with Engquist-Osher flux and min-mod-limiter (cf. \cite{Kroner1997}) to compute $\bar\varphi^\ksemi$ as the solution of
  \begin{equation}\label{eq:eoflux}
 \abs{Z_i} \frac{\bar\varphi^\ksemi_i - \bar\varphi^k_i}{\tau} + \sum_{j \in \mathcal{N}_i} F_{ij}(\vflow, \varphi^k_i, \varphi^k_j) = 0
 \end{equation}
 where $F_{ij}(\vflow, \cdot, \cdot)$ denotes the Engquist-Osher-flux \cite[3.2.6]{Kroner1997} which is Lip\-schitz continuous, antisymmetric and which satisfies the consistency condition 
\begin{equation*}
 F_{ij}(\vflow, \varphi, \varphi) = \skla{\nu_{ij}, \vflow} \varphi.
\end{equation*}
  Now, we transform $\bar\varphi^\ksemi$ back to obtain a corresponding finite element function $\varphi^\ksemi \in \Uh$. This way, $\varphi^\ksemi$ is an approximate solution of the convective part
  \begin{equation*}
 \frac{\varphi^\ksemi_i - \varphi^k_i}{\tau} + \Ph \skla{\vflow, \nabla \varphi^k} = 0.
\end{equation*}
 Then, we compute $\varphi^{k+1} \in \Uh$ as the solution of the purely diffusive Cahn-Hilliard problem 
\begin{subequations}
   \begin{equation*}
 \iOmega \rkla{\varphi^\kimp-\varphi^\ksemi} \psi + \tau\iOmega M(\varphi^k) \skla{\nabla\mu^\kimp, \nabla\psi} = 0 \quad \text{ for all }\psi \in \Uh
\end{equation*}
\begin{multline*}
 \iOmega \mu^\kimp \psi = \iOmega \skla{\nabla \varphi^\kimp, \nabla\psi} \\
  + \iOmega \Jh \rkla{\rkla{F'_+(\varphi^\kimp) + F'_-(\varphi^k)} \psi} \quad \text{ for all }\psi \in \Uh.
\end{multline*}
\end{subequations}
To handle the nonlinearity in $F'_+$, we use Newton's method and proceed along standard strategies for Cahn-Hilliard or lubrication type equations (see \cite{Grun2003} and the references therein).

\paragraph{{\bf Problem (ME)}}\label{par:problem_me}\addcontentsline{toc}{subsubsection}{Problem (ME)}
Let us compute new iterates $\vflow^{k+1}$ and $p^\kimp$ for velocity and pressure under the assumption that we are given a velocity iterate $\vflow^k$ corresponding to the previous time-step as well as phase-field iterates $\varphi^k$ and $\varphi^{k+1}$ corresponding to the old and the new time-step. 
Recall that this way we also have $\rho^k := \rho(\varphi^k)$ and $\rho^{k+1} := \rho(\varphi^{k+1})$ at hand.

We determine $\vflow^{k+1}$ as the solution of \eqref{eq:schemeB_a} which is obviously linear in $\vflow^{k+1}$. The first two terms on the left-hand side can be further simplified and rewritten as follows
\begin{equation}\label{eq:bk1}
I+II=\frac{1}{2\tau} \int_\Omega \rkla{\rho^k+\rho^\kimp} \Jhh \skla{\vflow^\kimp,\wflow} - \frac1\tau \iOmega \rho^k\Jhh \skla{\vflow^k,\wflow}.
\end{equation}
The second term contributes to the right-hand side of the linear equation for $\vflow^{k+1}$. 
Specifying a nodal basis of $\Uh$ as $\set{\psi_1,\dots, \psi_{\dim \Uh}}$ and a nodal basis of $\Xflow_h$ as $\set{\wflow_1,\dots, \wflow_{\dim \Xflow_h}}$, we expand
\begin{equation*}
\vflow^{k+1} := \sum^{\dim \Xflow_h}_{i=1} V_i^{k+1} \wflow_i \quad \text{and} \quad p^{k+1} := \sum^{\dim \Uh}_{i=1} P_i^{k+1} \psi_i.
\end{equation*}
Here, $V^{k+1}$ and $P^{k+1}$ are the coefficient vectors in $\R^{\dim \Xflow_h}$ and $\R^{\dim \Uh}$, respectively. Hence, we may rewrite the combination of \eqref{eq:schemeB_a} and \eqref{eq:schemeB_aa} as a saddle-point problem
\begin{multline} \label{eq:saddle}
  \begin{pmatrix}
   \frac{1}{2\tau}(\Mflow(\rho^k)+\Mflow(\rho^\kimp)) + \Aflow(\varphi^k) +\Nflow_a + \Nflow_b & \Bflow^T\\
    \Bflow & 0 \\
 \end{pmatrix}\begin{pmatrix} V^\kimp \\ P^\kimp
  \end{pmatrix} = \\ 
  \begin{pmatrix} \frac{1}{\tau}\Mflow(\rho^k) V^k +
    K \\ 0 
  \end{pmatrix}.
\end{multline}
Here, the matrix $\Mflow(\rho)$ is given by
\begin{equation*}
 \rkla{\Mflow(\rho)}_{ij}:=\iOmega\rho\Jhh \skla{\wflow_i,\wflow_j}
\end{equation*}
and $\Nflow_a=\Nflow_a(\rho, \vflow)$ is defined by
\begin{equation*} \rkla{\Nflow_a(\rho,\vflow)}_{ij}:=-\tfrac12\iOmega\rho \skla{\vflow,\rkla{\nabla\wflow_i}^T\wflow_j}
+ \tfrac12\iOmega\rho \skla{\vflow, \rkla{\nabla\wflow_j}^T\wflow_i}.
\end{equation*}  

Similarly, $\Nflow_b=\Nflow_b(\jflow)$ reflects the sum of the fifth and sixth term on the left-hand side in \eqref{eq:schemeB_a}, $\Aflow(\varphi)$ corresponds to the penultimate term on the left-hand side, and $K \in \R^{\dim \Xflow_h}$ is given by
\begin{equation*} 
 K_i = \rkla{K(\mu, \varphi, \kflow_\text{grav}) }_i := \iOmega \mu \skla{\nabla\varphi, \wflow_i} + \iOmega\skla{\kflow_\text{grav}, \wflow_i}.
 \end{equation*} 
Finally, $\Bflow \in \R^{\dim \Uh \times \dim \Xflow_h}$ is given by $\Bflow_{ij} := \iOmega \skla{\wflow_j, \nabla \psi_i}$ and corresponds to the solenoidality condition.

Numerical integration is exact for elementwise linear or bilinear functions and their products. Nonlinear coefficients are approximated by interpolation: $\eta(\varphi) \approx \Ih \eta(\varphi) = \sum_{i=0}^{\dim \Uh} \eta \rkla{\varphi(x_i)} \psi_i$.

We solve the saddle-point problem \eqref{eq:saddle} by Picard iteration (see \cite[7.2.2]{Elman2005}), which is in fact a fixed-point iteration in $V^{k+1}$ and $P^{k+1}$. 

Concerning discretization with respect to time, \eqref{eq:schemeB_a} indicates already our use of an semi-implicit Euler scheme.

Note also that the standard Taylor-Hood element (with polynomial degree of two and one) may be replaced by a stabilized P1-P1-approach first mentioned in \cite{Dohrmann2004}. Although inf-sup-stability is not shown, the results for equal polynomial degree in velocity and pressure are satisfactory, and computation time is tremendously decreased.

Now, all the informations are available to present the details of our splitting scheme.

\paragraph{{\bf Splitting Algorithm}}\label{par:splitting_algorithm}
Assume $(\vflow^k, \varphi^k, p^k)$ to be given in the $k$-th time-step. Define tolerances $\varepsilon_\vflow, \varepsilon_\varphi>0$.

\begin{enumerate}
 \item Set $i:=1$. Solve Problem (CH) using $\vflow^k$ and $\varphi^k$ to get an inner iterate $\varphi^{k+1, i}$.
 \item
  \begin{enumerate}
   \item Given $\varphi^{k+1, i}$ and $\vflow^k$, solve Problem (ME) to get $\vflow^{k+1, i}$ and $p^{k+1,i}$.
   \item Solve Problem (CH) using $\vflow^{k+1, i}$ and $\varphi^{k+1, i}$ to get $\varphi^{k+1, i+1}$. 
   \item If $\abs{\vflow^{k+1, i+1} - \vflow^{k+1, i}} > \varepsilon_\vflow$ or  $\abs{\varphi^{k+1, i+1} - \varphi^{k+1, i}} > \varepsilon_\varphi$, set $i:=i+1$ and go to (a).
  \end{enumerate}
  \item Take $\vflow^{k+1} := \vflow^{k+1, i+1}$, $p^{k+1} := p^{k+1, i+1}$ and $\varphi^{k+1} := \varphi^{k+1, i+1}$.
\end{enumerate}

\subsection{Adaptivity}\label{subsec:adaptivity}
 To enhance the performance of the scheme, we make use of adaptivity concepts both in space and time.
 Let us begin with adaptivity in space. We refine the grid by bisection. Since no a posteriori-error estimators are available, we use the moduli of $\nabla \varphi$, $\nabla \vflow_1$ and $\nabla \vflow_2$ as control parameters. This way, we may guarantee that the diffuse interface, i.e. the narrow strip around the zero-level line of $\varphi$, is resolved by sufficiently many grid points. In practice, we aim at around 20 grid points perpendicular to the level line.
 
 To this scope, we proceed as follows. We compute $M := \max\limits_{K \in \Th}\abs{\nabla \varphi_{|K}}$ and $m := \min\limits_{K \in \Th}\abs{\nabla \varphi_{|K}}$. Such elements $K$ for which
\begin{equation*} \abs{\nabla \varphi_{|K}} > (1-C^\text{ref}) \cdot m + C^\text{ref} \cdot M, \quad \quad C^\text{ref}=0.1,
  \end{equation*}
 are marked for refinement. Elements $K$ for which 
\begin{equation*}
 \abs{\nabla \varphi_{|K}} < (1-C^\text{coarse}) \cdot m + C^\text{coarse} \cdot M, \quad \quad C^\text{coarse}=0.2,
   \end{equation*}
 are marked for coarsening.
 
 Similar control mechanisms apply to $\nabla \vflow_{1}$ and $\nabla \vflow_{2}$. However, the constants are different, we take $C^\text{ref} = 0.1$, $C^\text{coarse}=0.5$. Note that we always give preference to refinement.

 Concerning adaptivity with respect to time, we follow the idea in \cite{Grun2000} selecting the time-increment to be controlled by the ratio
  \begin{equation*}\frac{\text{grid size}}{\text{propagation speed of the interface}}.
    \end{equation*}
  In our case, we take $\nabla \mu$ and $\vflow$ as estimators for the propagation speed. We determine the increment $\tau_k$ in the $k$-th time-step by
  \begin{equation}\label{eq:timeincrement}
   \tau_k := 0.9 \cdot h \cdot \rkla{\max\set{ \min\set{ \max\limits_{K \in \Th}\set{\abs{\nabla\mu^k_{|K}}, \abs{\vflow^k_{|K}}}, \vflow_\text{max}}, \vflow_\text{min}}}^{-1}.
   \end{equation}
 We use cut-off values $\vflow_\text{min} = 10$ and $\vflow_\text{max} = 10^5$ to bound $\tau_k$ by positive constants.

\subsection{Coding}\label{subsec:coding}
 We implement the algorithm in C++, using the Intel Math Kernel Library\footnote{Intel MKL: \url{software.intel.com/en-us/articles/intel-mkl/}}. It provides the PARDISO solver (\cite{Schenk2004}) which is used to solve the linearized saddle-point problems. We prefer it compared to UMFPACK or Krylov-space-methods for reasons of performance. For the linearized 4-th order problem, we use BiCGstab. For visualization, we use ParaView\footnote{ParaView: \url{http://www.paraview.org/}}.

\section{Numerical examples}\label{sec:results}
In this section, we validate the numerical scheme presented before. The goal is three-fold.

First, we give strong indication for the convergence of the scheme when the gridsize is decreasing (see Subsection \ref{subsec:eoc}).
Secondly, we study the influences different density ratios (expressed by the Atwood number $A= \frac{\hat\rho_1 - \hat\rho_2}{\hat\rho_1 +
  \hat\rho_2}$) have on the qualitative behaviour of solutions as well as on their stability (see Subsection \ref{subsec:rising_droplet}).
 Moreover, we present a simulation on Rayleigh-Taylor instability (see Subsection \ref{subsec:rayleigh-taylor}).
Finally, we give an example indicating that there may be drastic differences between solutions computed for model \cite{DSS07} and solutions computed for model \cite{AGG}.

We mention that $\delta=0.05$  and $M=0.005$ are taken in all the experiments presented in this paper, if not specified differently.

\subsection{Experimental order of convergence}\label{subsec:eoc}
\begin{figure}[ht]
  \centering
  \psfrag{rxline}{${{r_x}}$}
  \psfrag{ryline}{$r_y$}
  \psfrag{time}{time}
  \psfrag{majorminorradiusplaceplace}{major / minor halfaxes}
  \includegraphics[width=0.6\textwidth,angle=270]{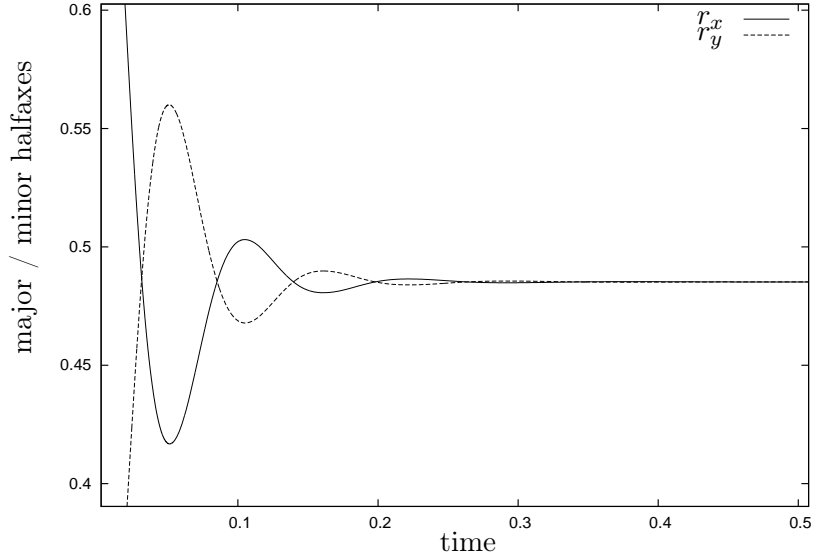}
 \caption{ Retraction of an ellipsoidal droplet to a circular shape. Radii $r_x$ and $r_y$  (level 14, Taylor-Hood elements)}
  \label{fig:retraction}
 \end{figure}

 We consider two liquids with equal viscosities ($\eta=0.01$) and densities $\hat\rho_1=0.001$, $\hat\rho_2=0.019$ (which correspond to an Atwood number of 0.9). We choose a constant mobility $M=0.5$ and the interface parameter $\delta=0.1$. We assume external forces to vanish, and as initial data we take an ellipsoidal droplet of liquid 2 with halfaxes $r_x \approx 0.87$ and $r_y \approx 0.29$ located in the center of the quadratic domain $\rkla{-1,1}^2$ and surrounded by liquid 1.
 Initially, we assume the configuration to be at rest, and we watch the droplet deforming to a circular shape. Due to inertia, it overshoots but reaches a stable circular shape after some time (numerically stationary after $t>1.5$), see Fig.~\ref{fig:retraction}.
 
As a reference solution, we use computations on a uniform spatial grid, where time increments are chosen according to \eqref{eq:timeincrement}. For stabilized P1-P1 elements, we allow for refinement level 16 (Table~\ref{table:eoc_p1}), for Taylor-Hood elements we take level 14 (Table~\ref{table:eoc_p2}). We compare solutions computed on uniform grids of different refinement levels and a solution computed using adaptive grid refinement with these reference solutions. For the errors of the phase-field measured in the $L^2$-norm, we find a decay which hints at an experimental order of convergence above 2.

\begin{table}[ht]
  \begin{center}
    \begin{tabular}{c c|| l l l}
      Level&h& T=0.4 & T=1.0 & T=2.5 \\
      \hline
      10 &0.0625& $5.03750e^{-2}$ & $5.03860e^{-2}$ & $5.03865e^{-2}$ \\
      12 &0.03125& $7.49162e^{-3}$ & $7.42532e^{-3}$ & $7.42515e^{-3}$ \\
      14 &0.015625& $1.56562e^{-3}$ & $1.98117^{e-3}$ & $1.98117e^{-3}$ \\
      \hline
      adaptive &from 0.0625 & \multirow{2}{*}{$2.04969e^{-3}$} & \multirow{2}{*}{$2.04688e^{-3}$} & \multirow{2}{*}{$2.02904e^{-3}$} \\
      (10-16)&to 0.0078125&\\
    \end{tabular}
        \caption{$L^2$-comparison against a reference solution (level 16, stabilized P1-P1 elements)}
  \end{center}
  \label{table:eoc_p1}
\end{table}

\begin{table}[ht]
  \begin{center}
    \begin{tabular}{c c|| l l l}
      Level&h & T=0.4 & T=1.0 & T=2.5 \\
      \hline
      10 &0.0625&$4.16917e^{-2}$ & $5.09377e^{-2}$ & $5.09385e^{-2}$ \\
      12 &0.03125& $8.11914e^{-3}$ & $7.98854e^{-3}$ & $7.98844e^{-3}$ \\
      \hline
      adaptive &from 0.0625 &\multirow{2}{*}{$5.00866e^{-3}$} & \multirow{2}{*}{$6.13341e^{-3}$} &\multirow{2}{*}{$1.43544e^{-2}$} \\
      (10-14) &to 0.015625 & \\
    \end{tabular}
        \caption{$L^2$-comparison against a reference solution (level 14, Taylor-Hood elements)}
  \end{center}
    \label{table:eoc_p2}
\end{table}

\subsection{Rising droplet -- Stability and tip formation}\label{subsec:rising_droplet}
In this series of experiments, we consider light droplets surrounded by an ambient liquid of lower viscosity and subjected to gravity forces. This choice of parameters becomes relevant in extraction processes of chemical engineering, for instance to extract a toluol or butanol droplet out of a water reservoir.

We identify the two-dimensional fluid domain with the set  $\Omega = [0,1] \times [0,2]$ in $\R^2$.  We assume it mainly to be occupied by fluid 1 ($\varphi \approx -1$) with the exception of a small circular-shaped subdomain at $(0.5, 0.5)$ of radius 0.5, where fluid 2 is found $(\varphi \approx +1)$. As values for the viscosity, we choose $\eta_1 = 0.001$, $\eta_2 \in \set{0.001, 0.1}$. We fix an average density $\hat\rho_\text{avg} =\frac12(\hat\rho_1 + \hat\rho_2)$ = 0.01 and let the  Atwood number $A= \frac{\hat\rho_1 - \hat\rho_2}{\hat\rho_1 + \hat\rho_2}$ range from $0.5$ to $0.99$, depending on the experiment. The force field is given by $\kflow_\text{grav} = (0,-10^4)^T$. At the liquid-solid interface, we have the choice between no-slip or Navier-slip boundary conditions. Since the domain $\Omega$ is comparatively small, we prefer Navier slip to minimize boundary effects. 

We use Taylor-Hood elements to investigate droplet shapes during the evolution depending on different Atwood numbers. In Fig.~\ref{fig:rising_e0001_at05} and \ref{fig:rising_e01_at05}, we choose the Atwood number to be 0.5. We observe droplet rising with slight deformation. It is not surprising that the deformation of the rising bubble is the stronger, the lower its viscosity liquid is chosen.

When higher Atwood numbers (0.9 in Fig.~\ref{fig:rising_e01_at09} and 0.99 in Fig.~\ref{fig:rising_e0001_at099}) are combined with small viscosities, we observe a pronounced tip-formation at the rising droplet. A liquid jet is found at the droplet front, the symmetry break-up apparently occurring in the experiments may be due to a flow which is no longer laminar around the tip.  

Let us emphasize that we use these last experiments to test the stability of the scheme in case of critical parameters. Presently, we are not aware of physical settings which would correspond to these choices (i.e. large discrepancy in mass densities \emph{and} the denser liquid being of lower viscosity).

 \newlength{\picscale}  

 \begin{figure}[ht]
  \centering
  \setlength{\picscale}{0.18\textwidth} 
  \fbox{\includegraphics[width=\picscale]{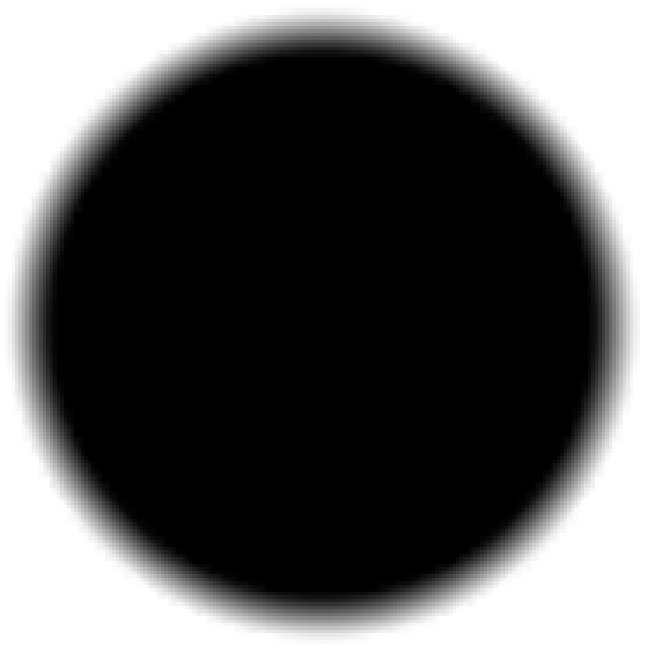}}
  \fbox{\includegraphics[width=\picscale]{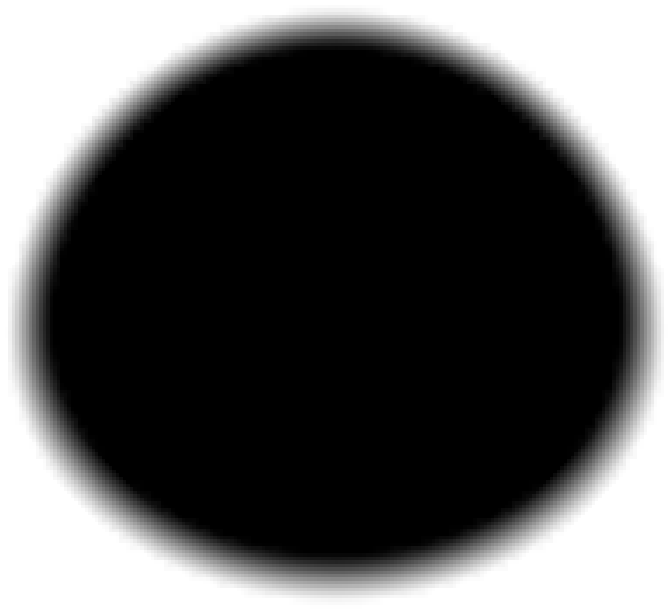}}
  \fbox{\includegraphics[width=\picscale]{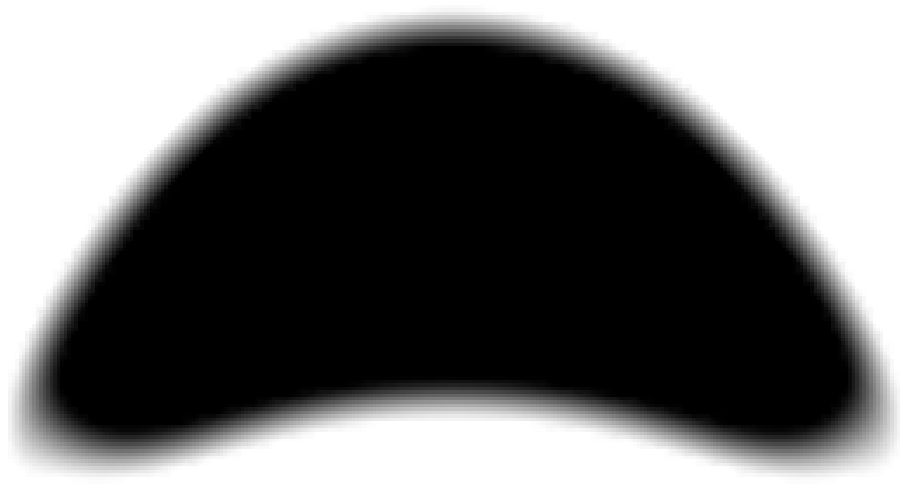}}
  \fbox{\includegraphics[width=\picscale]{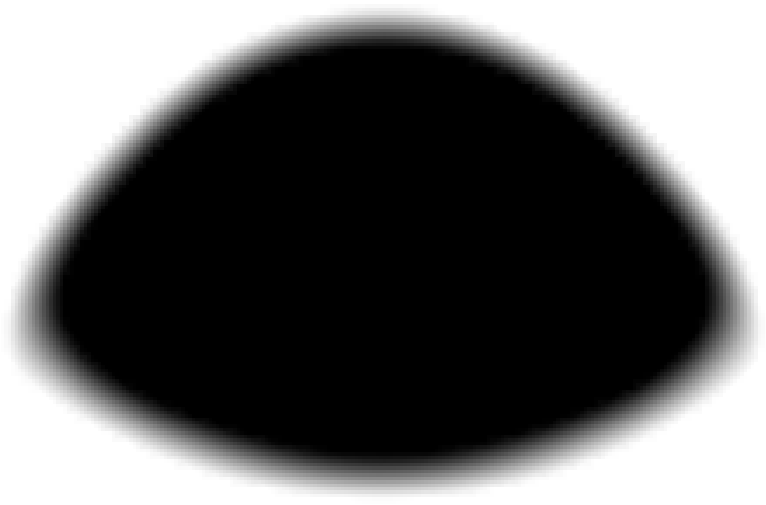}}
  \fbox{\includegraphics[width=\picscale]{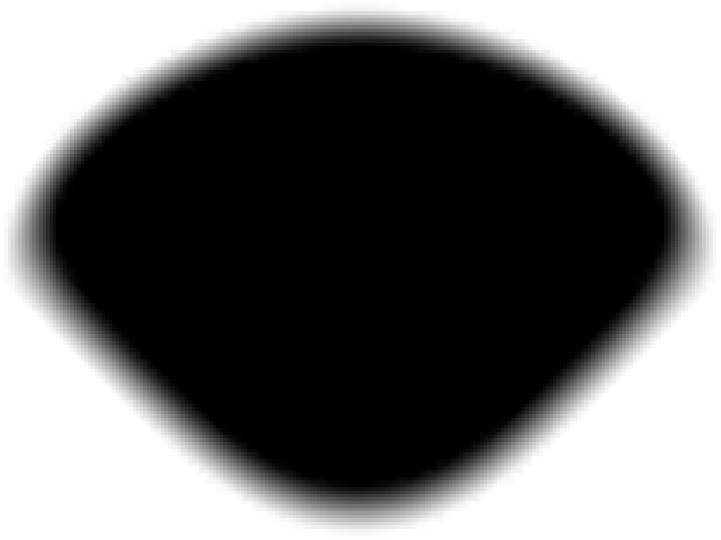}} 
  \fbox{\includegraphics[width=\picscale]{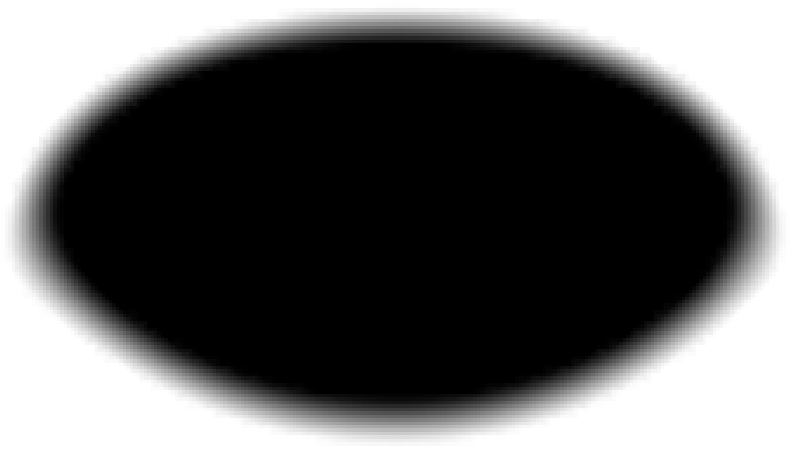}}
  \fbox{\includegraphics[width=\picscale]{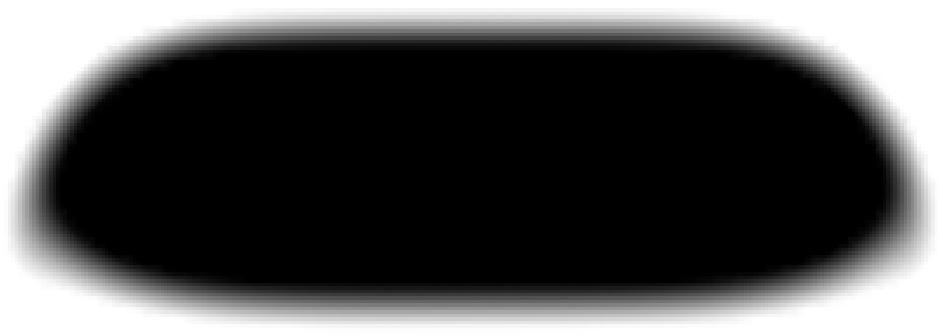}}
  \fbox{\includegraphics[width=\picscale]{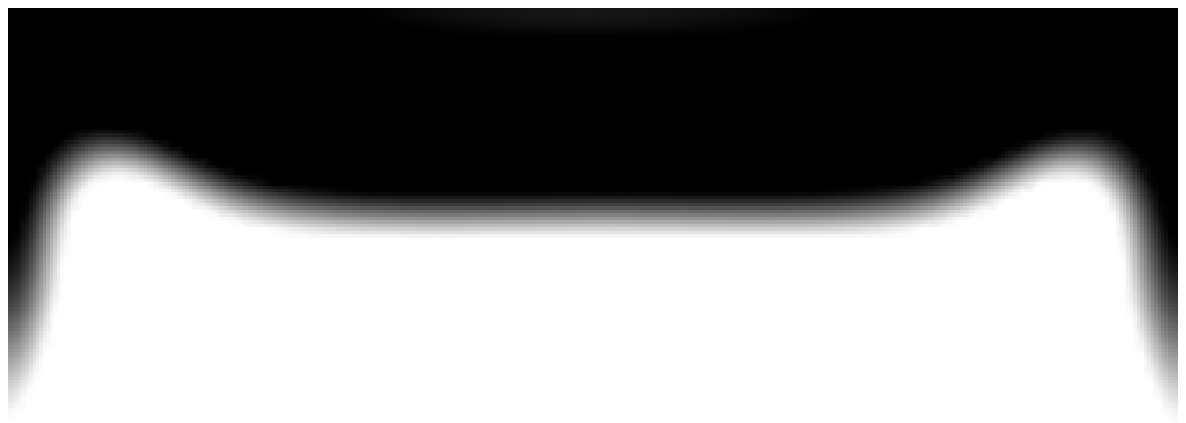}}
  \caption{Rising droplet in a constant gravitational field. Parameters $\hat\rho_\text{avg} = 0.01$, Atwood = 0.5, $\eta_1 = 0.001$, $\eta_2 = 0.001$ at time $T=$ 0, 0.01, 0.02, 0.03, 0.035, 0.04, 0.045, 0.05, $\delta=0.05$, $M =0.005$ (to be read linewise from top left to bottom right).}
  \label{fig:rising_e0001_at05}
 \end{figure}

 \begin{figure}[ht]
  \centering
  \setlength{\picscale}{0.18\textwidth} 
  \fbox{\includegraphics[width=\picscale]{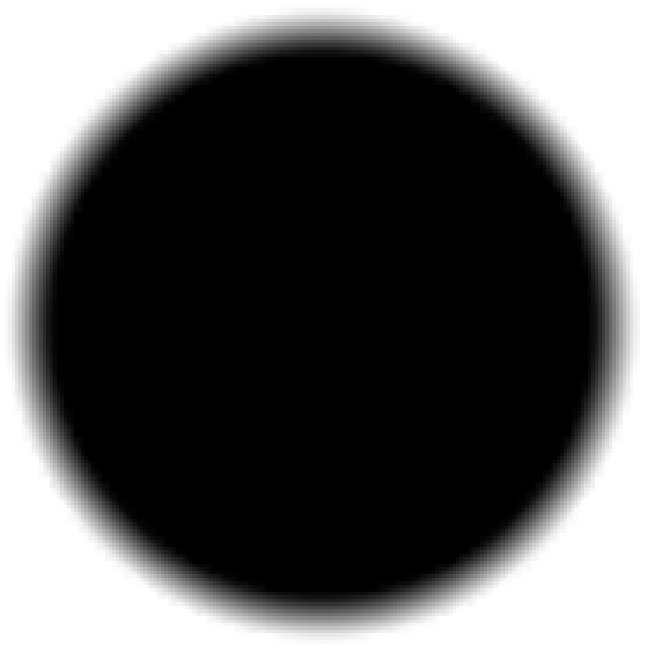}}
  \fbox{\includegraphics[width=\picscale]{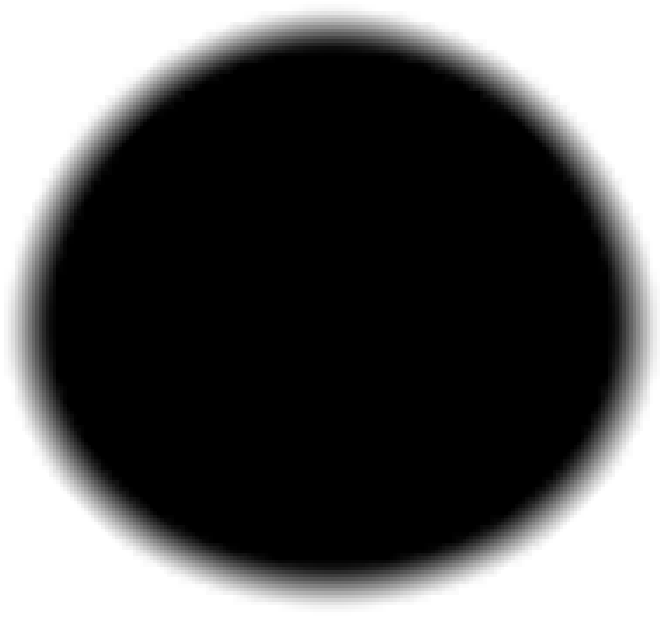}}
  \fbox{\includegraphics[width=\picscale]{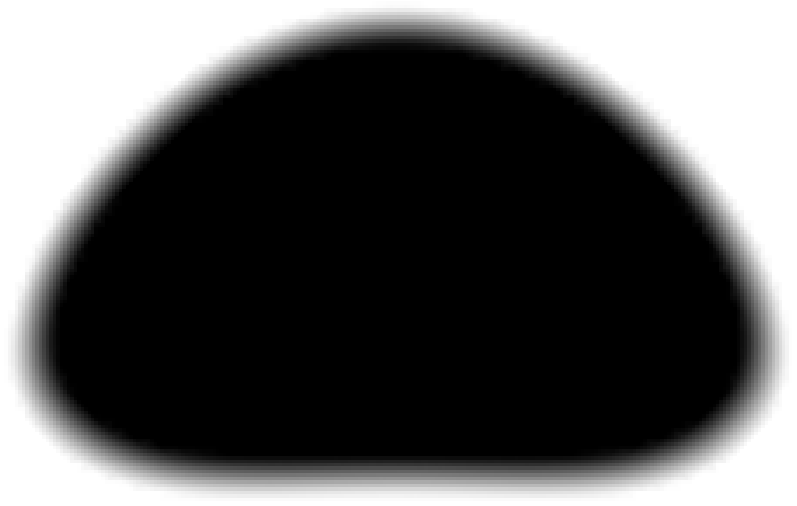}}
  \fbox{\includegraphics[width=\picscale]{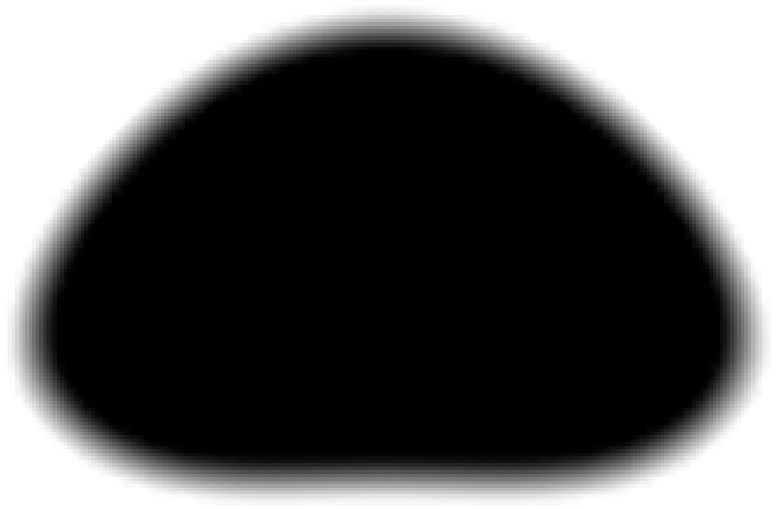}}
  \fbox{\includegraphics[width=\picscale]{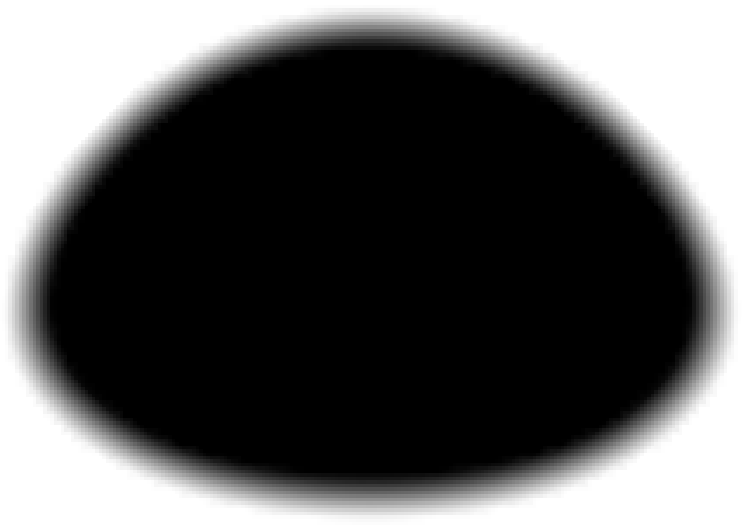}}  
  \fbox{\includegraphics[width=\picscale]{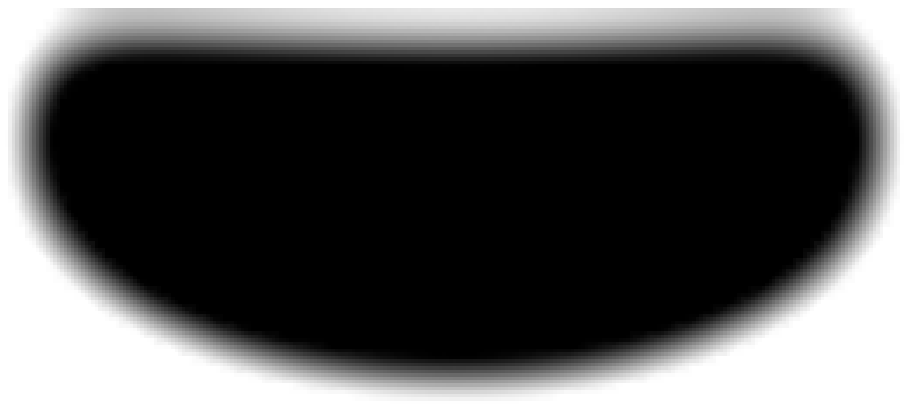}}
  \fbox{\includegraphics[width=\picscale]{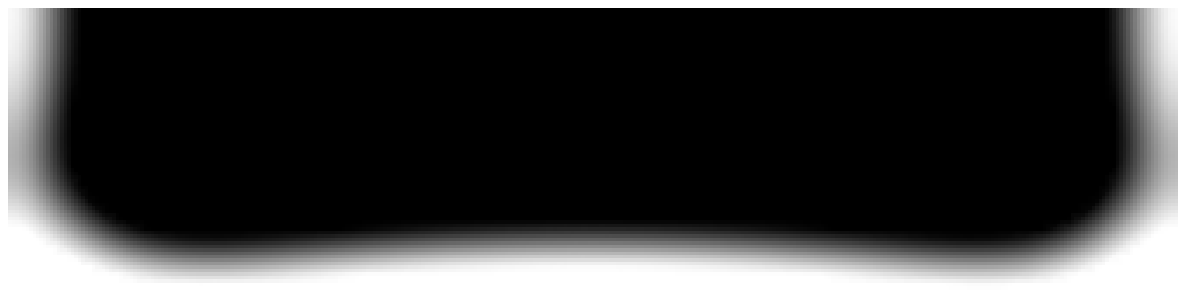}} 
  \fbox{\includegraphics[width=\picscale]{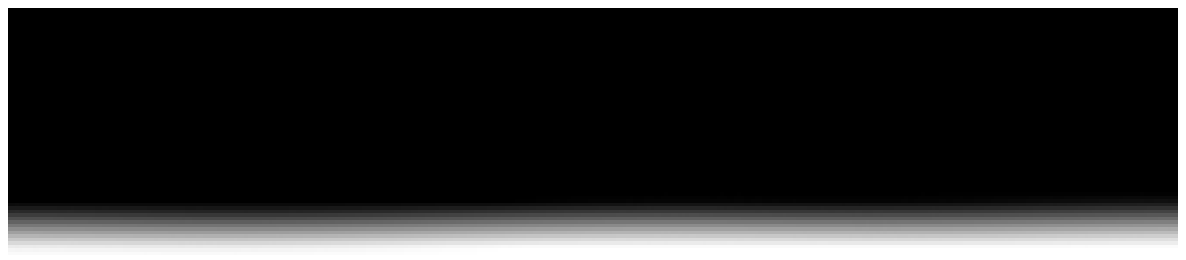}} 
  \caption{Rising droplet in a constant gravitational field. Parameters $\hat\rho_\text{avg} = 0.01$, Atwood = 0.5, $\eta_1 = 0.001$, $\eta_2 = 0.1$ at time $T=$ 0, 0.01, 0.02, 0.03, 0.04, 0.05, 0.06, 0.1, $\delta=0.05$, $M =0.005$ (to be read linewise).}
  \label{fig:rising_e01_at05}
 \end{figure}

 \begin{figure}[ht]
  \centering
  \setlength{\picscale}{0.14\textwidth}  
  \fbox{\includegraphics[width=\picscale]{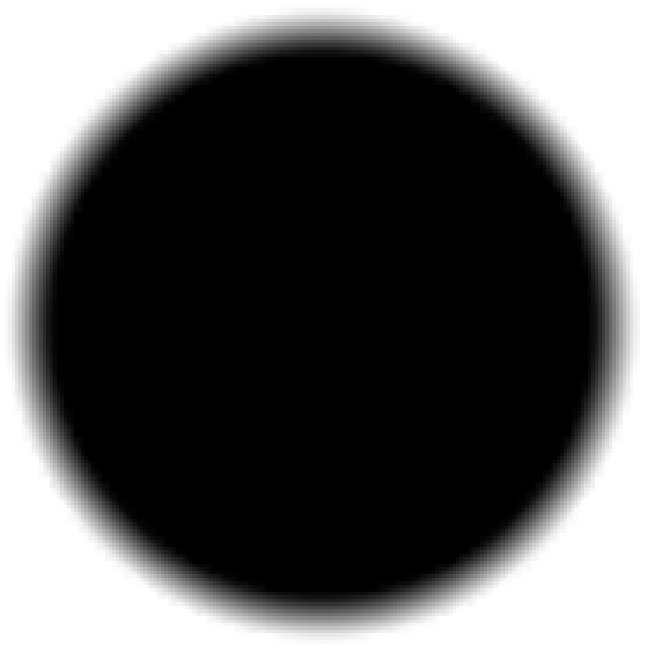}}
  \fbox{\includegraphics[width=\picscale]{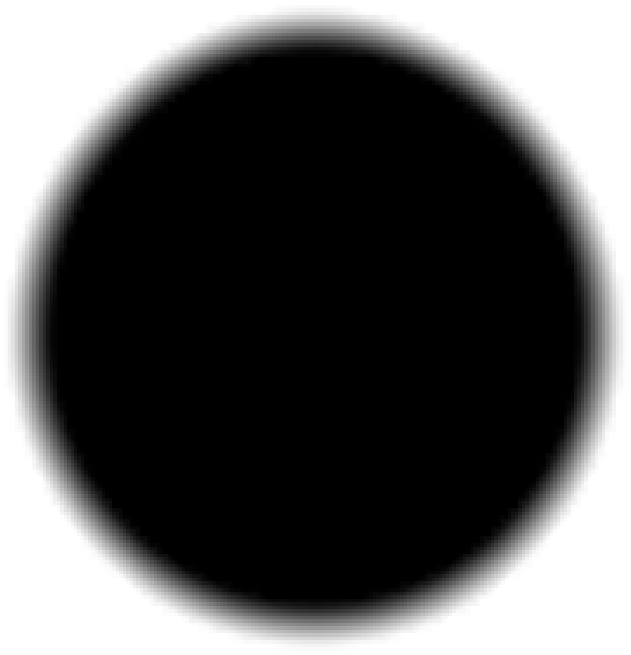}}
  \fbox{\includegraphics[width=\picscale]{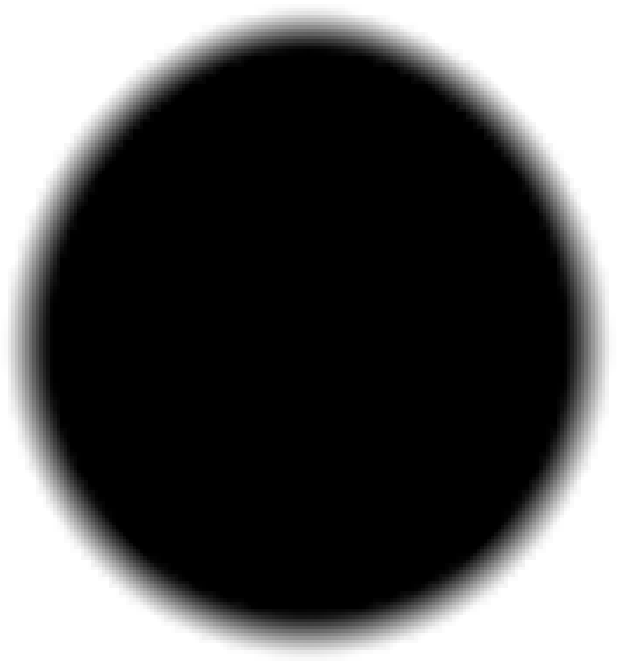}}
  \fbox{\includegraphics[width=\picscale]{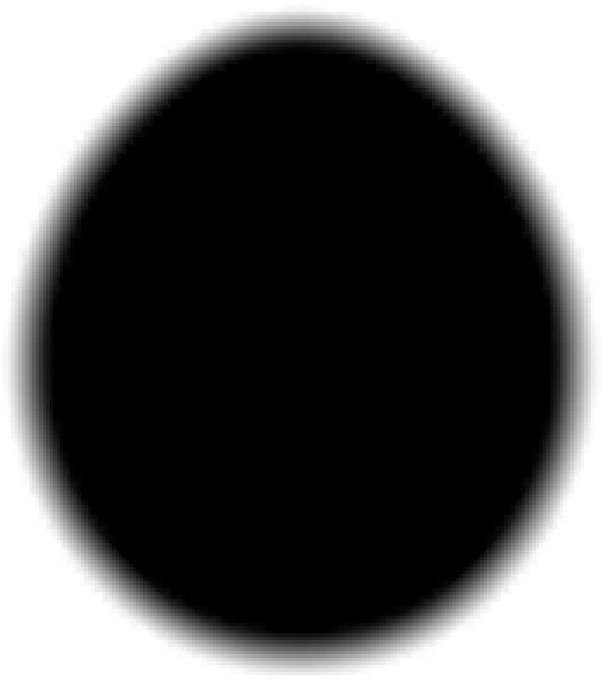}}
  \fbox{\includegraphics[width=\picscale]{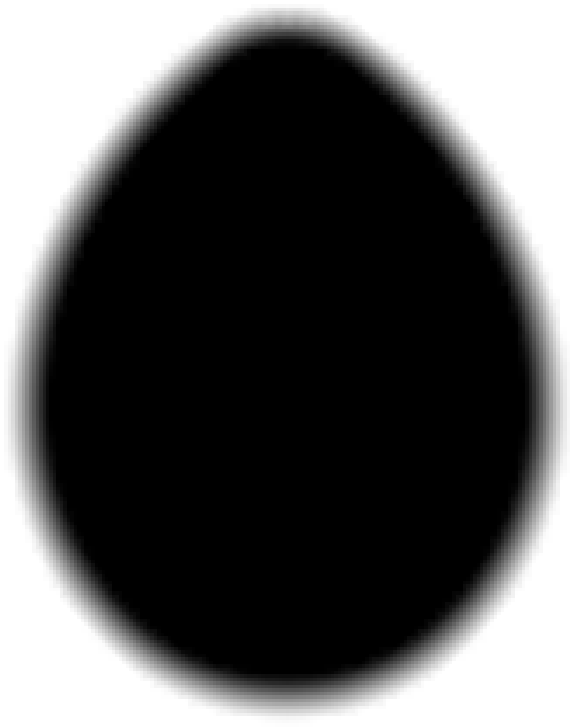}}
  \fbox{\includegraphics[width=\picscale]{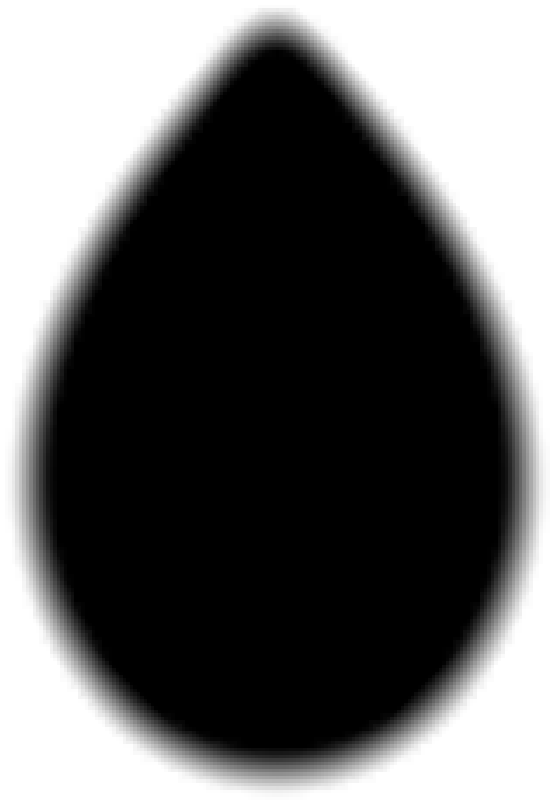}}
  \fbox{\includegraphics[width=\picscale]{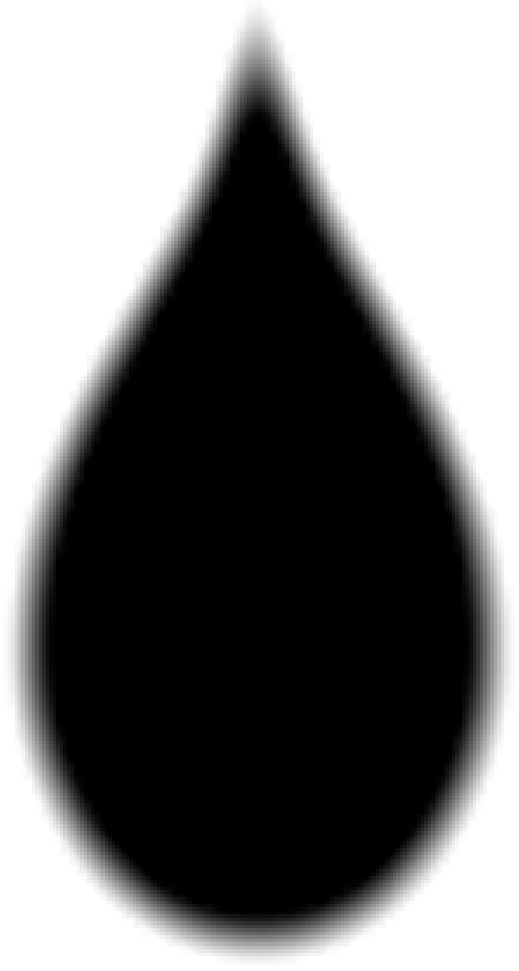}}
  \fbox{\includegraphics[width=\picscale]{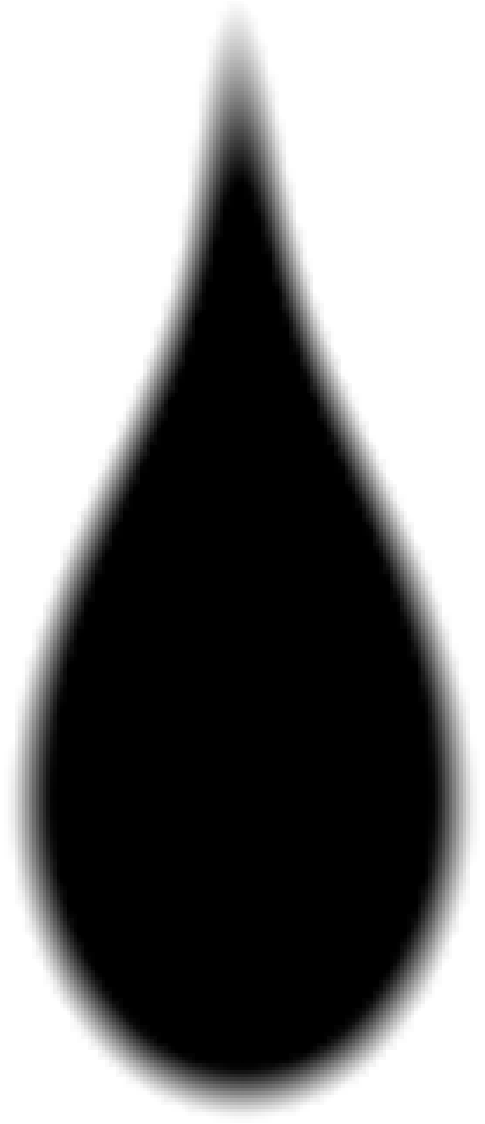}}
  \fbox{\includegraphics[width=\picscale]{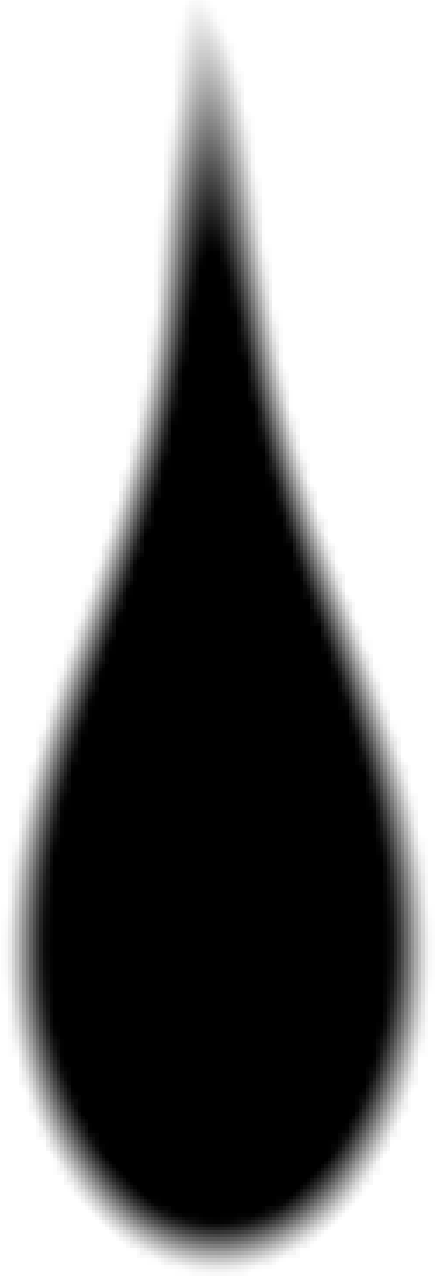}}
  \fbox{\includegraphics[width=\picscale]{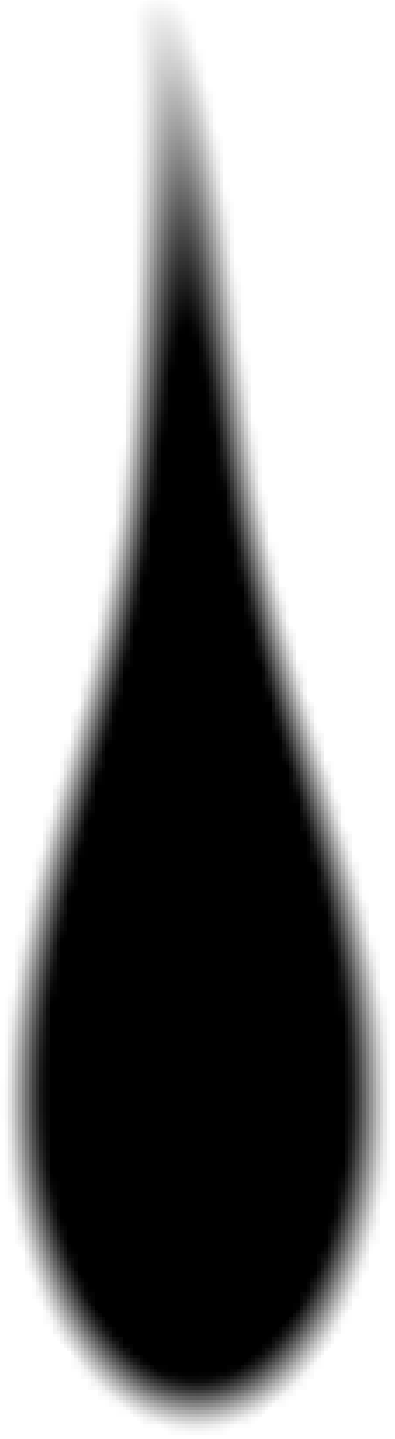}}
  \fbox{\includegraphics[width=\picscale]{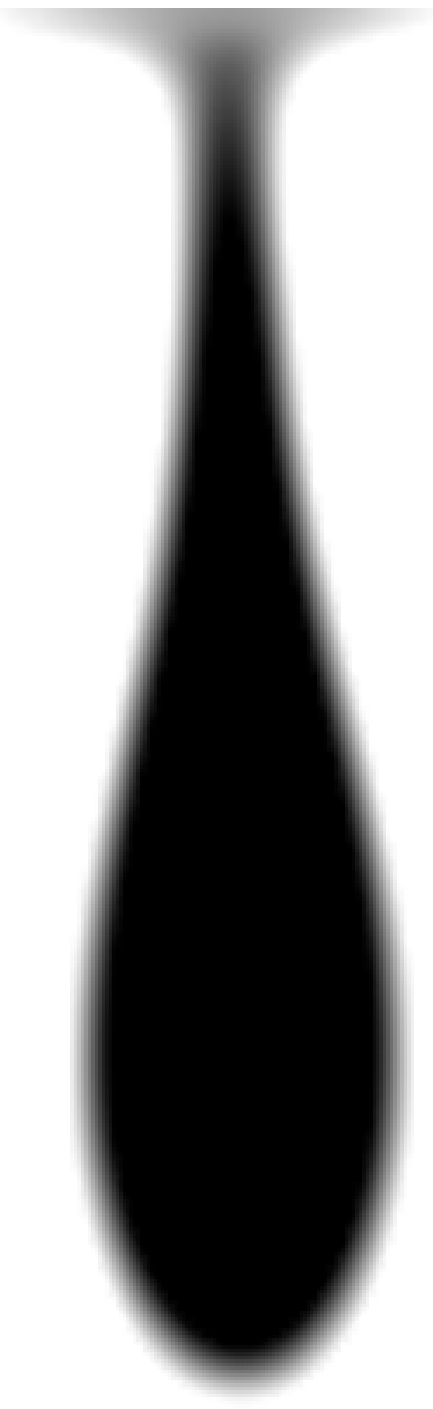}}
  \fbox{\includegraphics[width=\picscale]{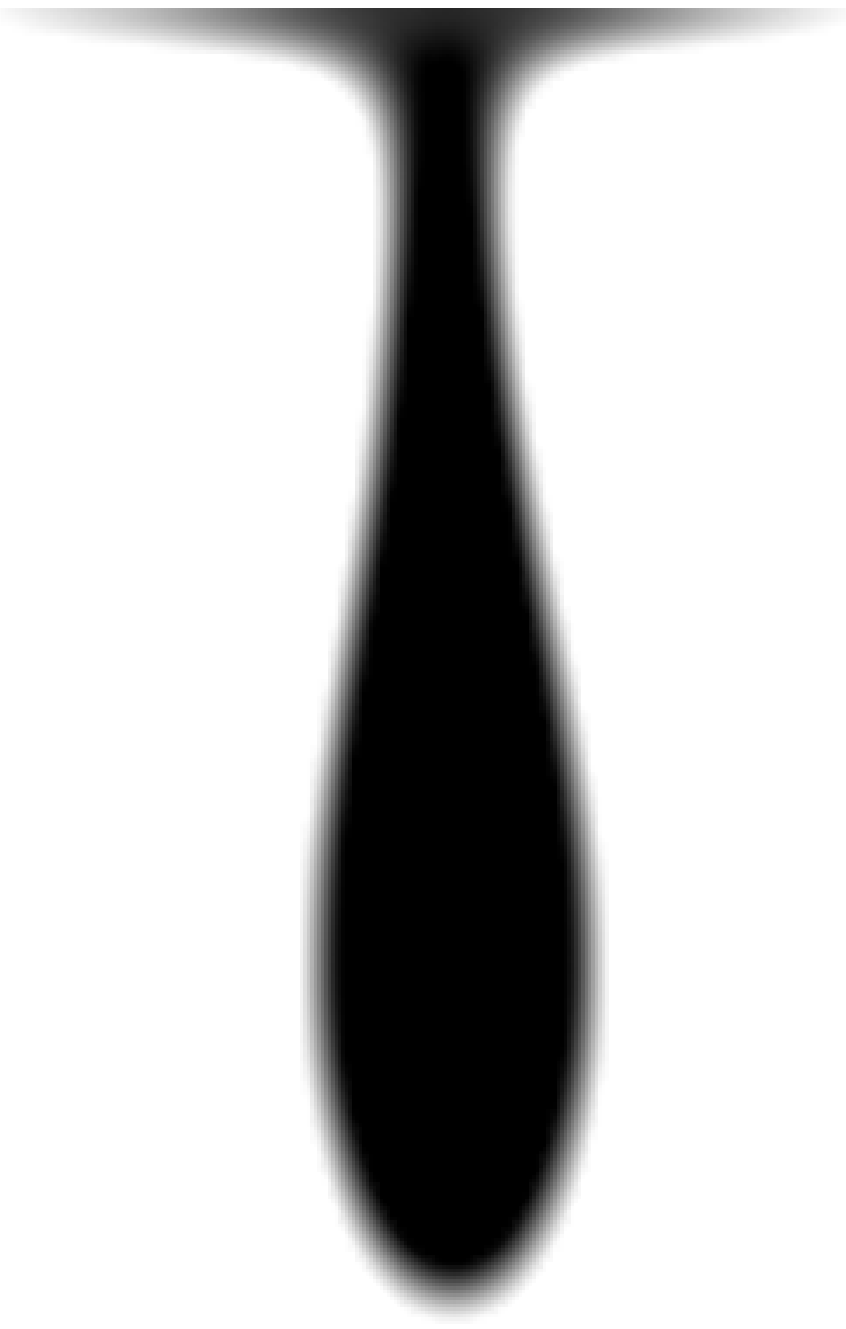}}
       
  \caption{Rising droplet in a constant gravitational field. Parameters $\hat\rho_\text{avg} = 0.01$, Atwood = 0.9, $\eta_1 = 0.001$, $\eta_2= 0.1$ at time $T=$ 0, 0.005, 0.006, 0.007, 0.008, 0.009, 0.01, 0.011, 0.012, 0.013, 0.014, 0.015, $\delta=0.05$, $M =0.005$ (to be read linewise).}
  \label{fig:rising_e01_at09}  
 \end{figure}

  \begin{figure}[ht]
  \centering
  \setlength{\picscale}{0.14\textwidth}
  \fbox{\includegraphics[width=\picscale]{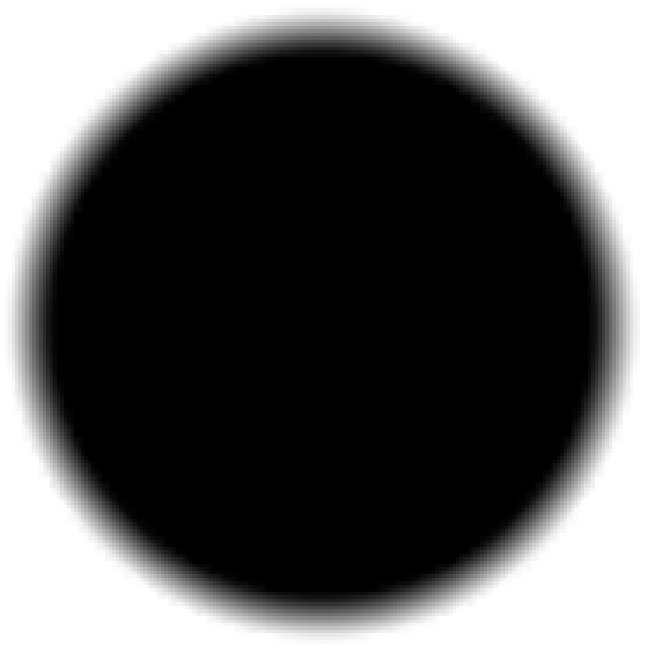}}
  \fbox{\includegraphics[width=\picscale]{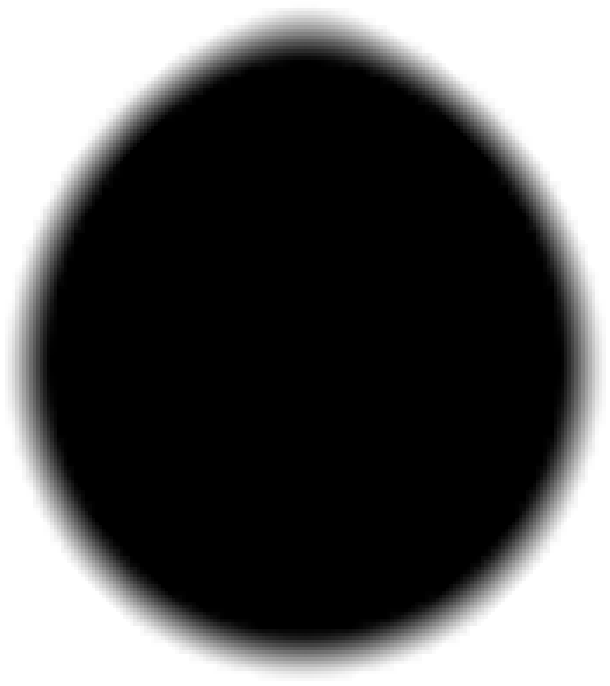}}
  \fbox{\includegraphics[width=\picscale]{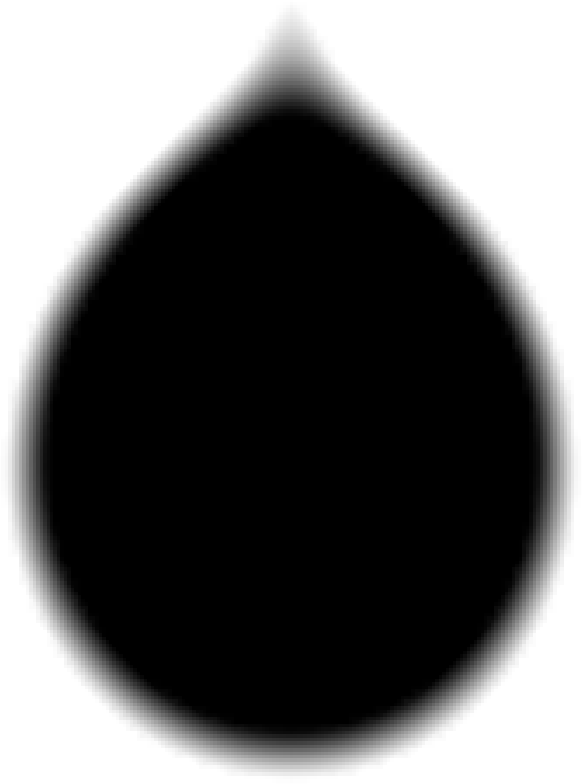}}
  \fbox{\includegraphics[width=\picscale]{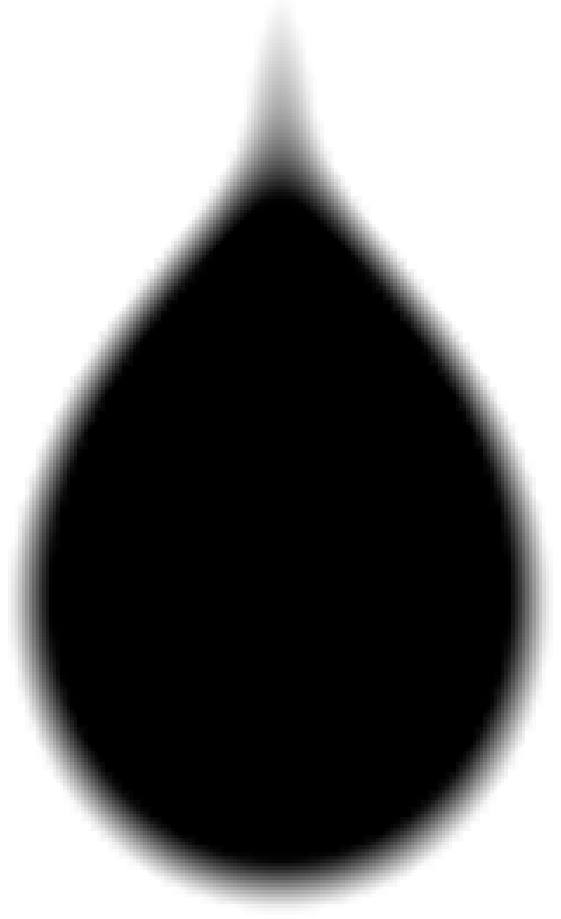}}
  \fbox{\includegraphics[width=\picscale]{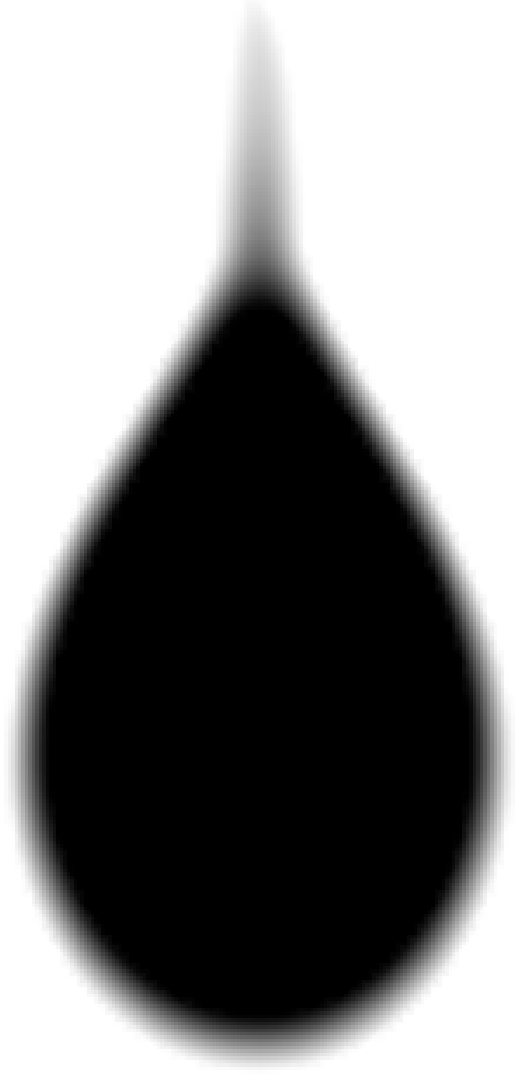}}
  \fbox{\includegraphics[width=\picscale]{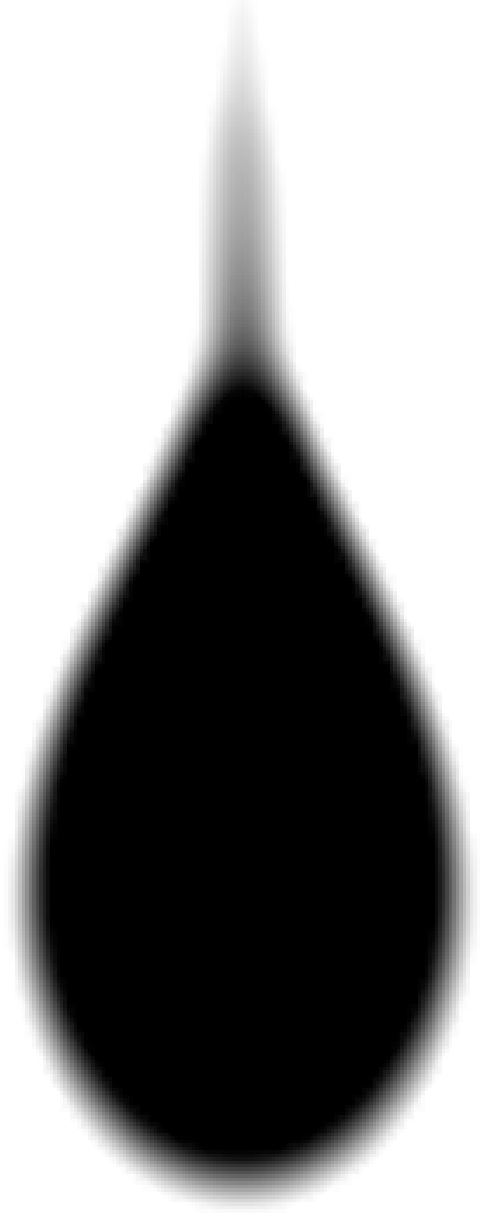}}
  \fbox{\includegraphics[width=\picscale]{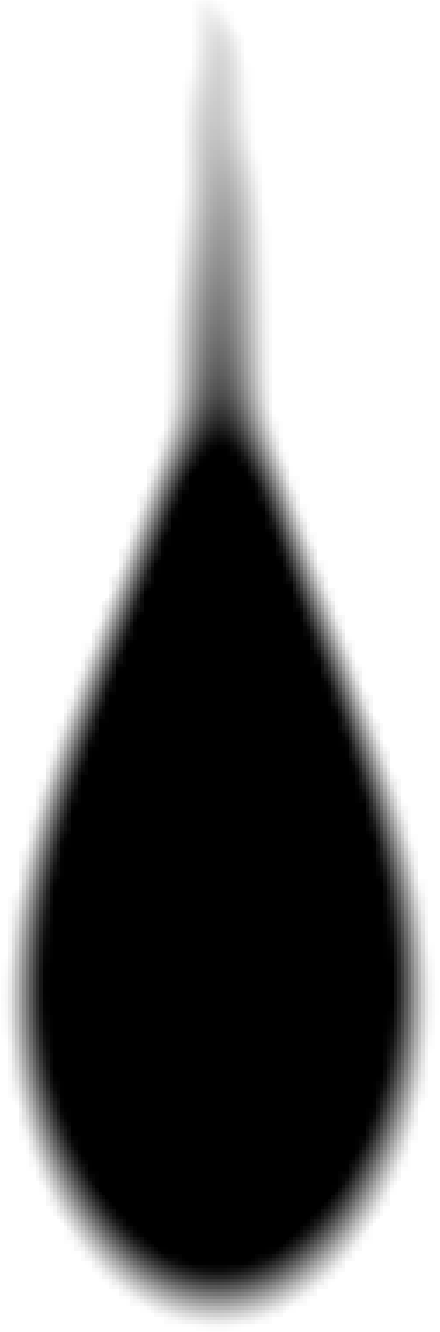}}
  \fbox{\includegraphics[width=\picscale]{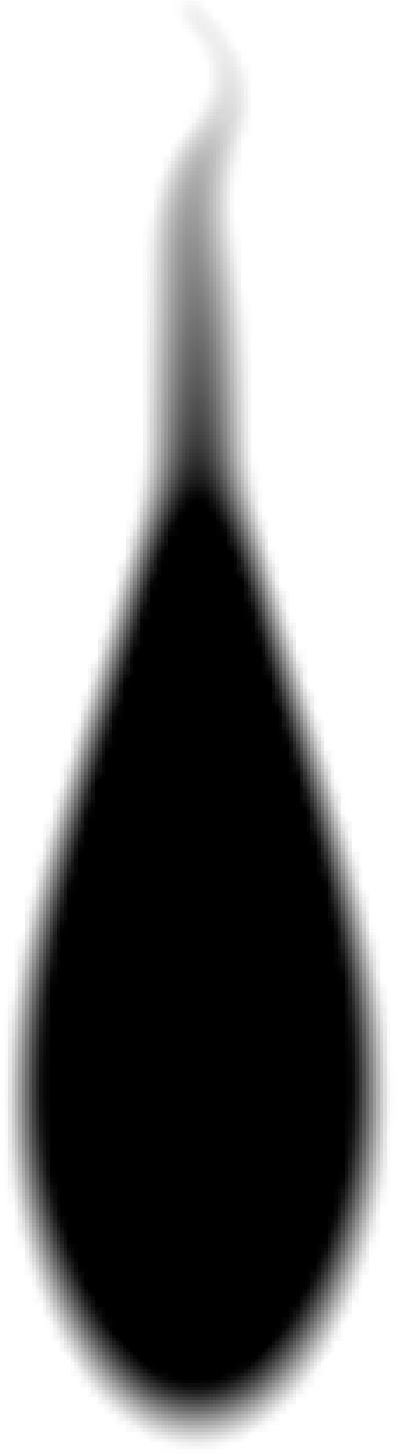}}
  \fbox{\includegraphics[width=\picscale]{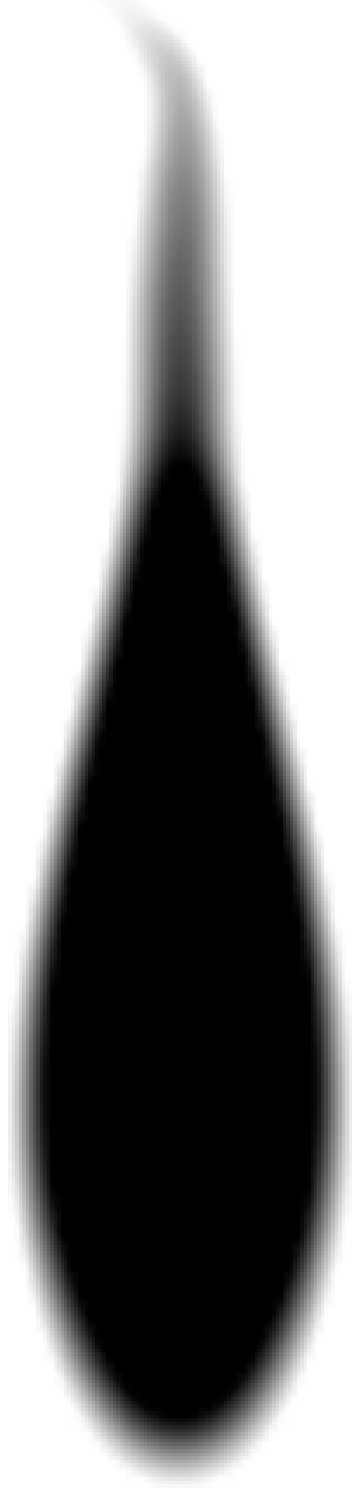}}
  \fbox{\includegraphics[width=\picscale]{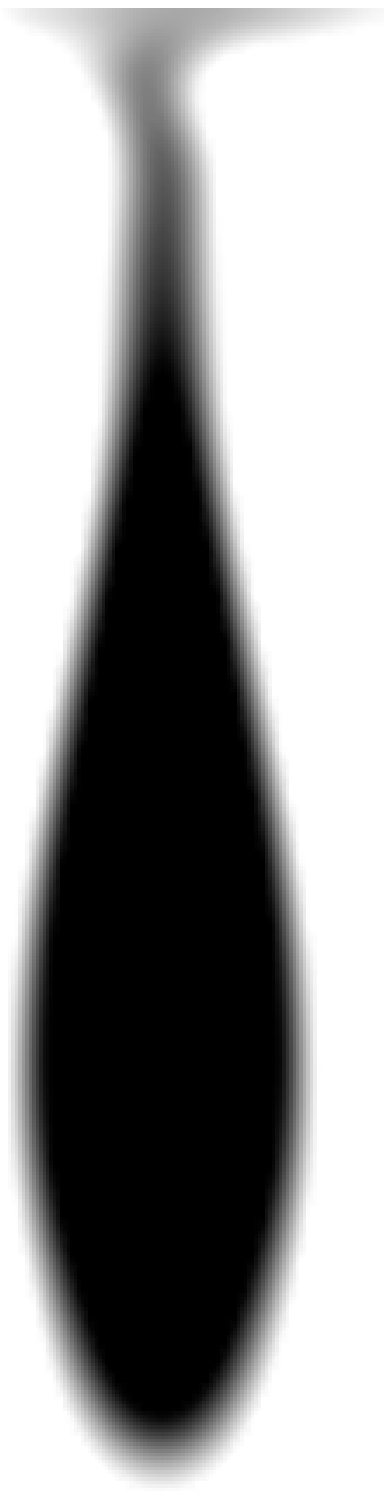}}
  \fbox{\includegraphics[width=\picscale]{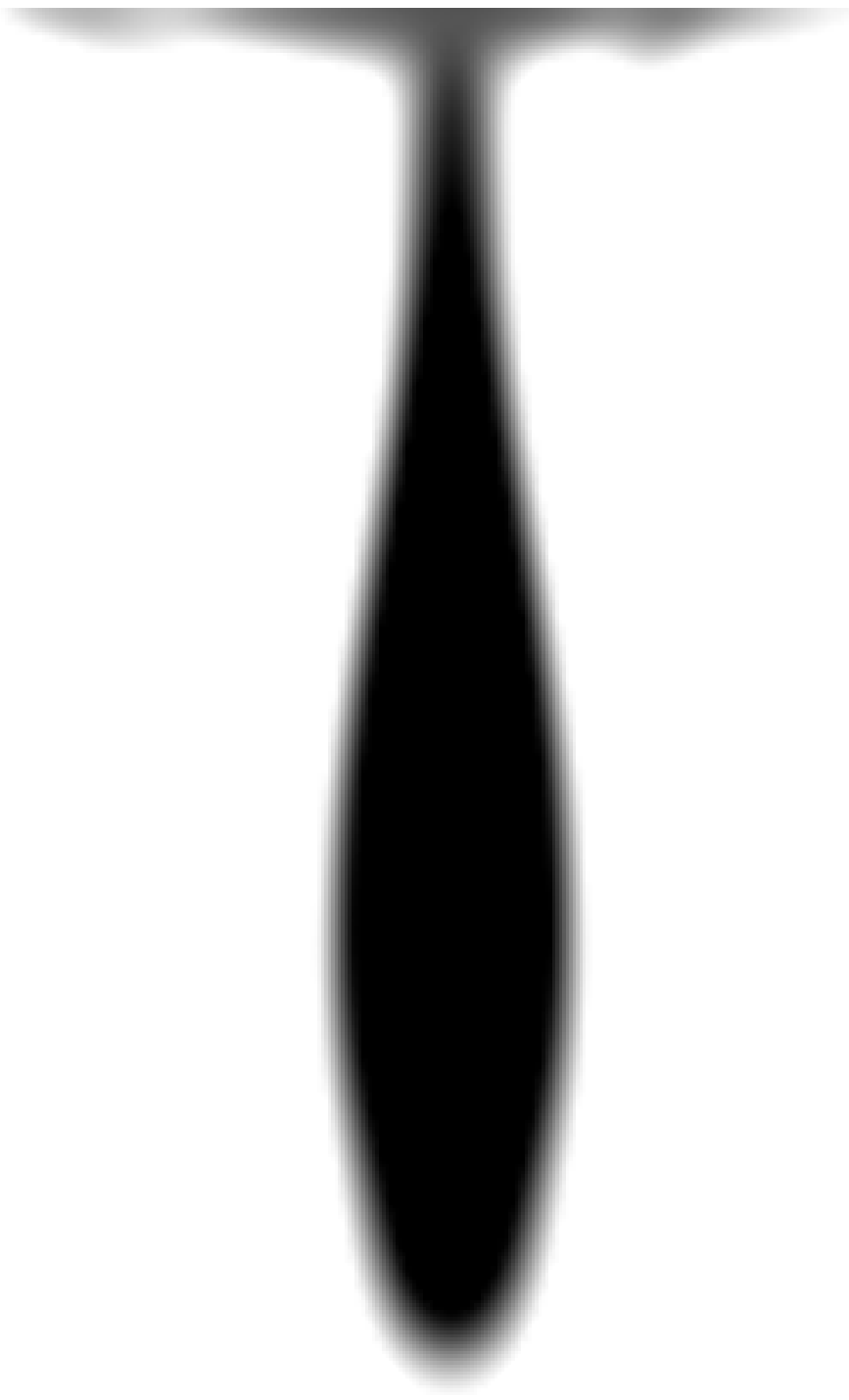}}     
  \fbox{\includegraphics[width=\picscale]{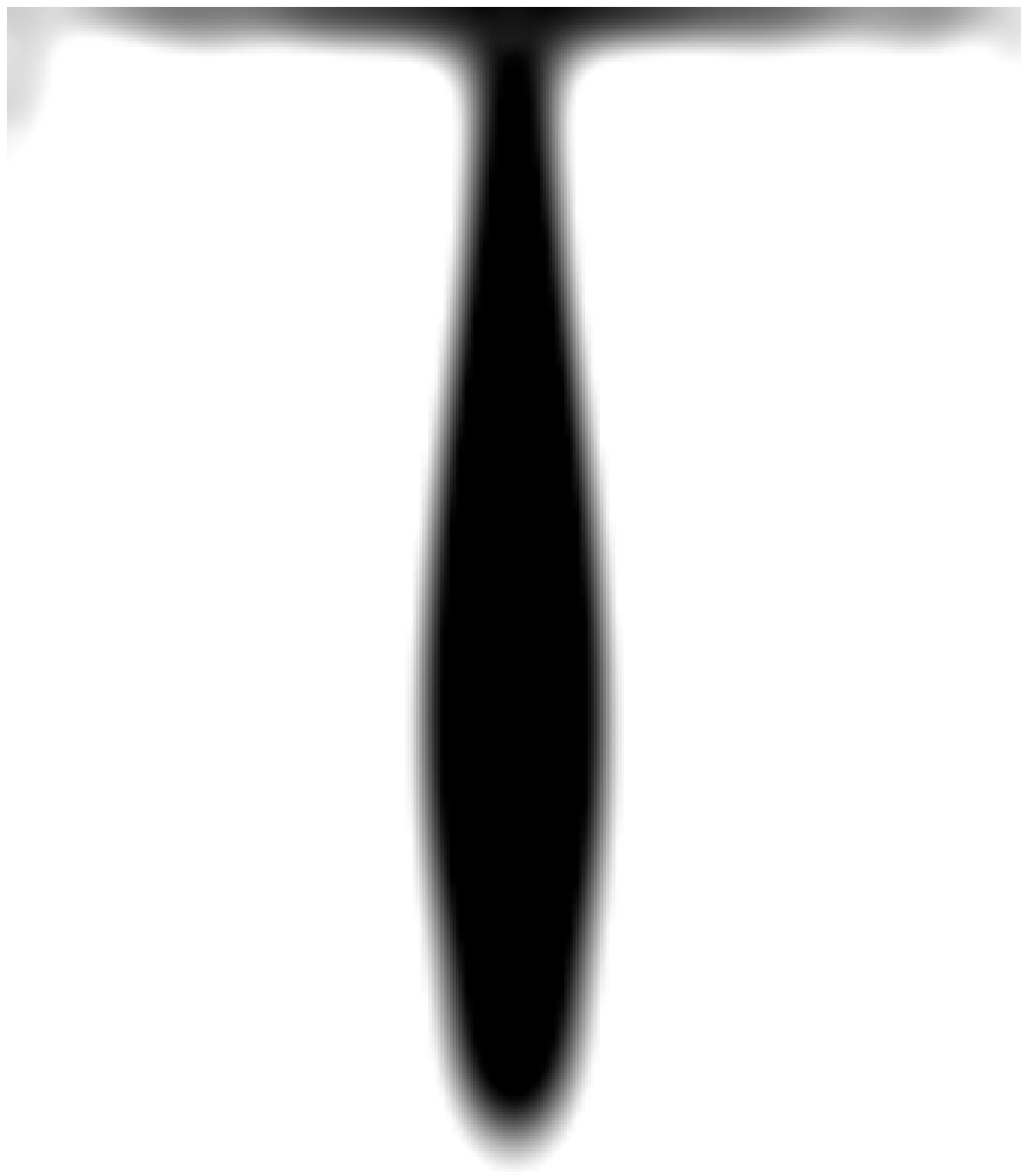}}  
  \caption{Rising droplet in a constant gravitational field. Parameters $\hat\rho_\text{avg} = 0.01$, Atwood = 0.99, $\eta_1 = 0.001$, $\eta_2 =0.001$ at time $T=$ 0, 0.005, 0.006, 0.007, 0.008, 0.009, 0.01, 0.011, 0.012, 0.013, 0.014, 0.015, $\delta=0.05$, $M =0.005$ (to be read linewise).}
  \label{fig:rising_e0001_at099}
 \end{figure}

\subsection{Rayleigh-Taylor instability}\label{subsec:rayleigh-taylor}
\begin{figure}[ht]
  \centering
 
  \setlength{\picscale}{0.25\textwidth}  
  \fbox{\includegraphics[width=\picscale]{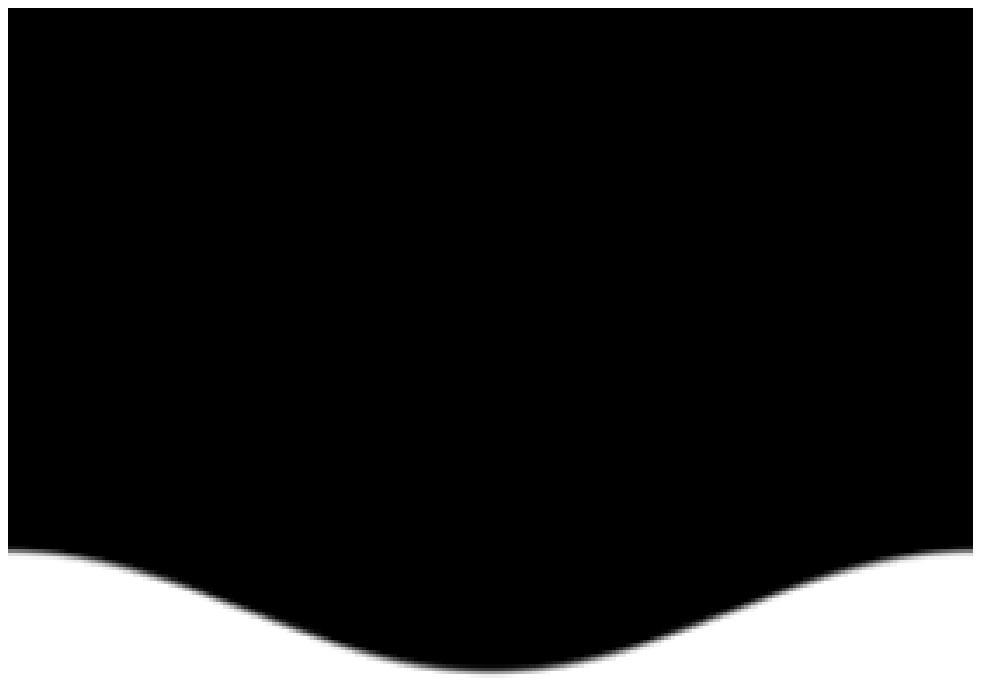}}  
  \fbox{\includegraphics[width=\picscale]{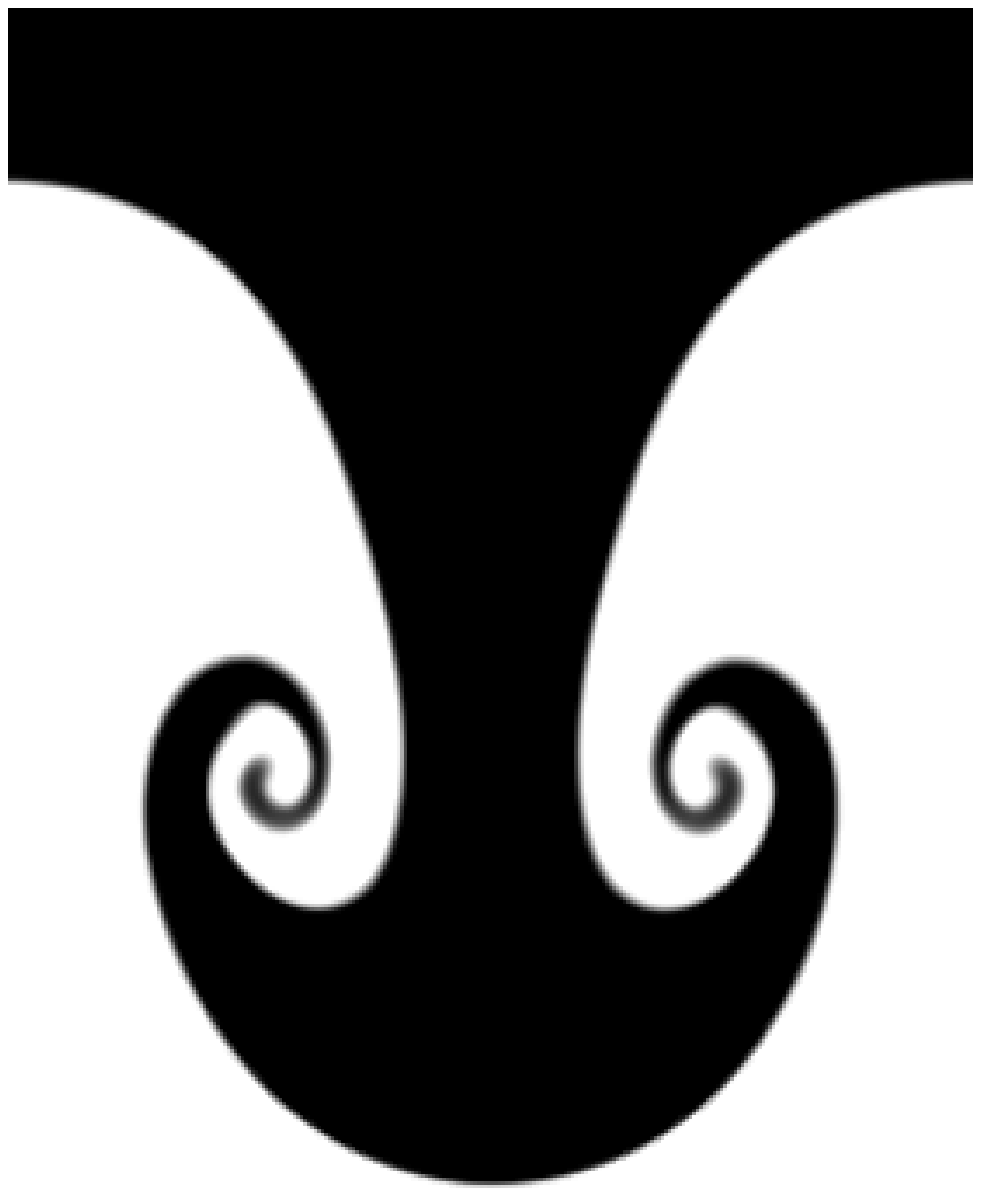}}
  \fbox{\includegraphics[width=\picscale]{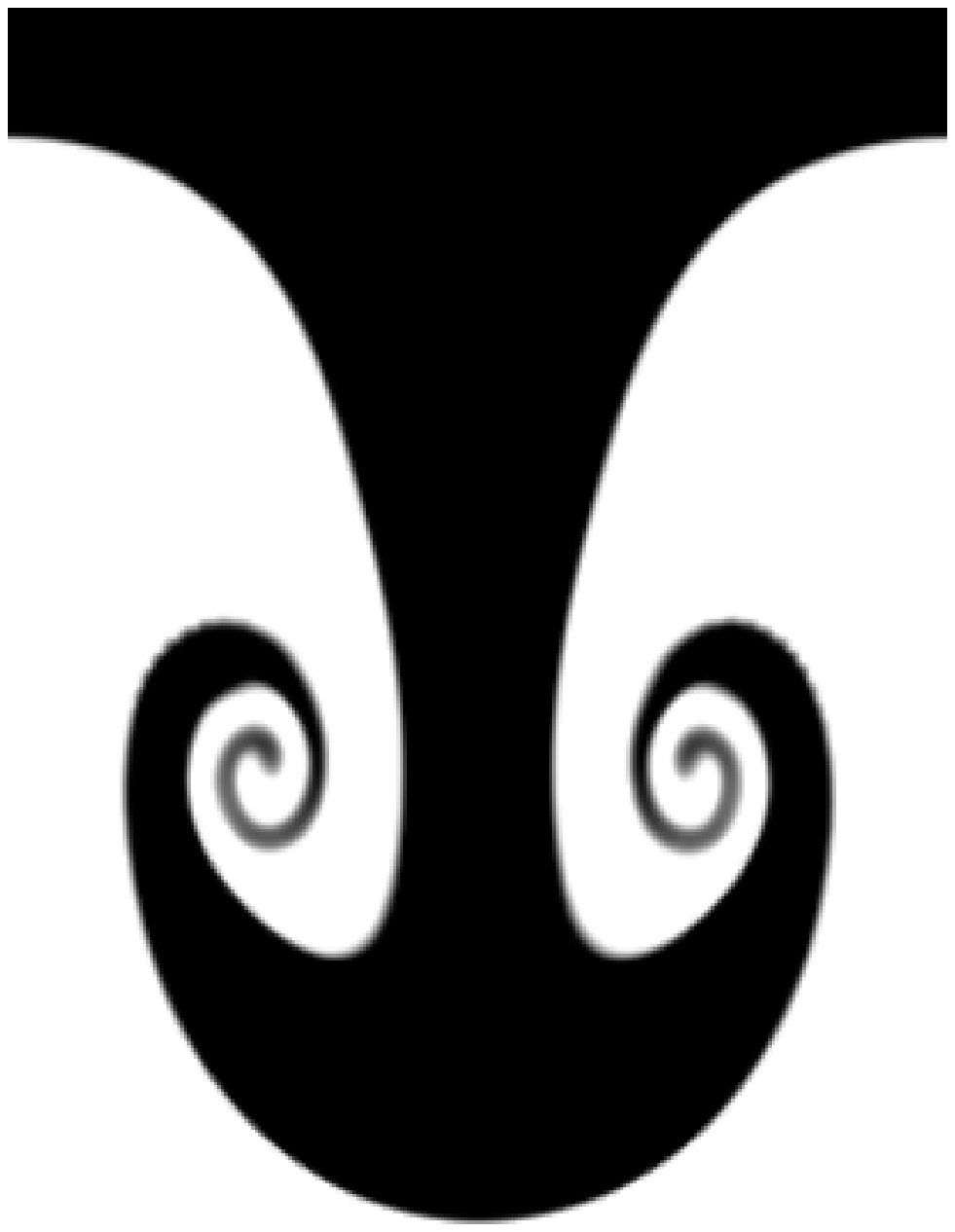}}
  \fbox{\includegraphics[width=\picscale]{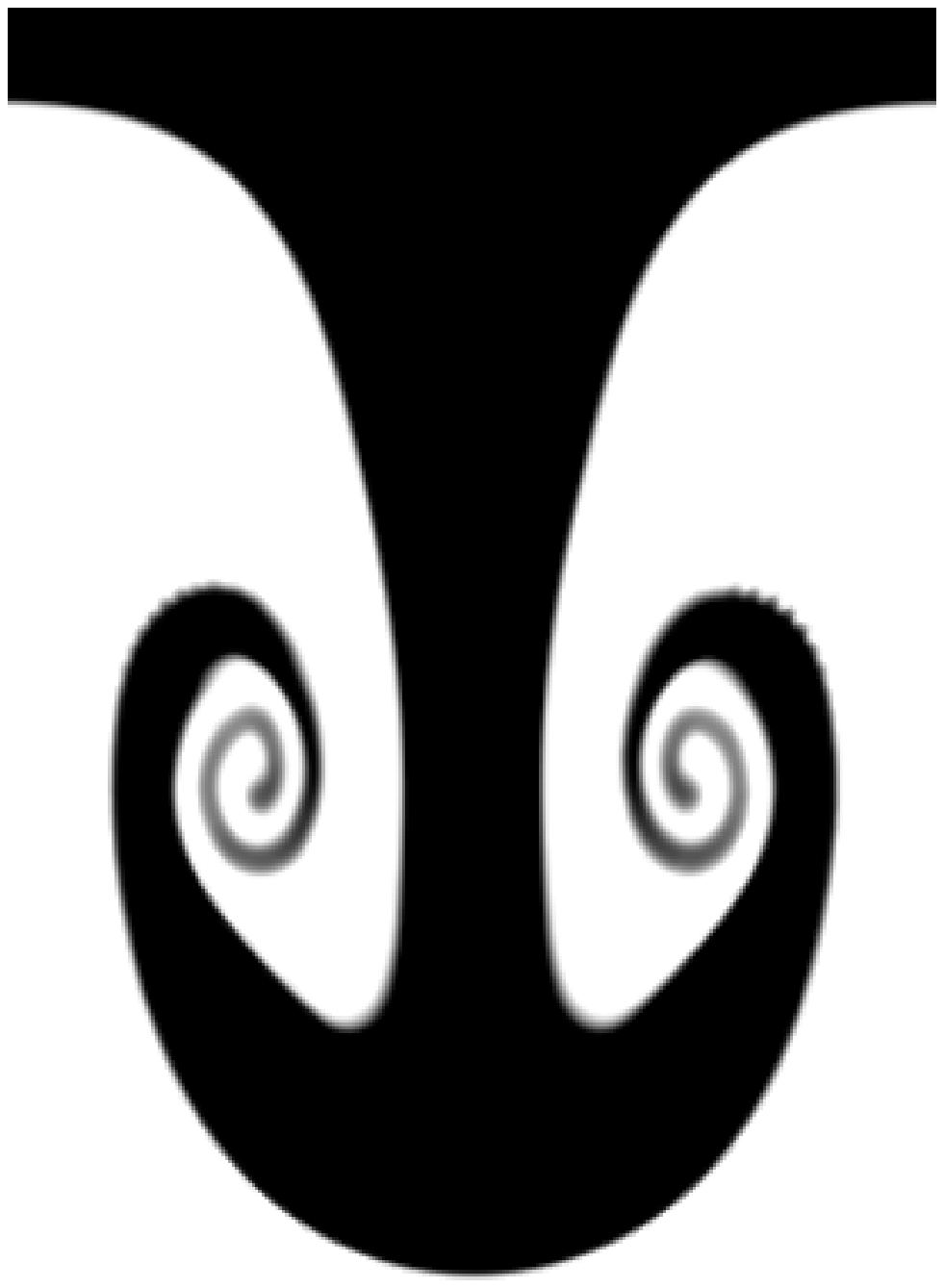}}
  \fbox{\includegraphics[width=\picscale]{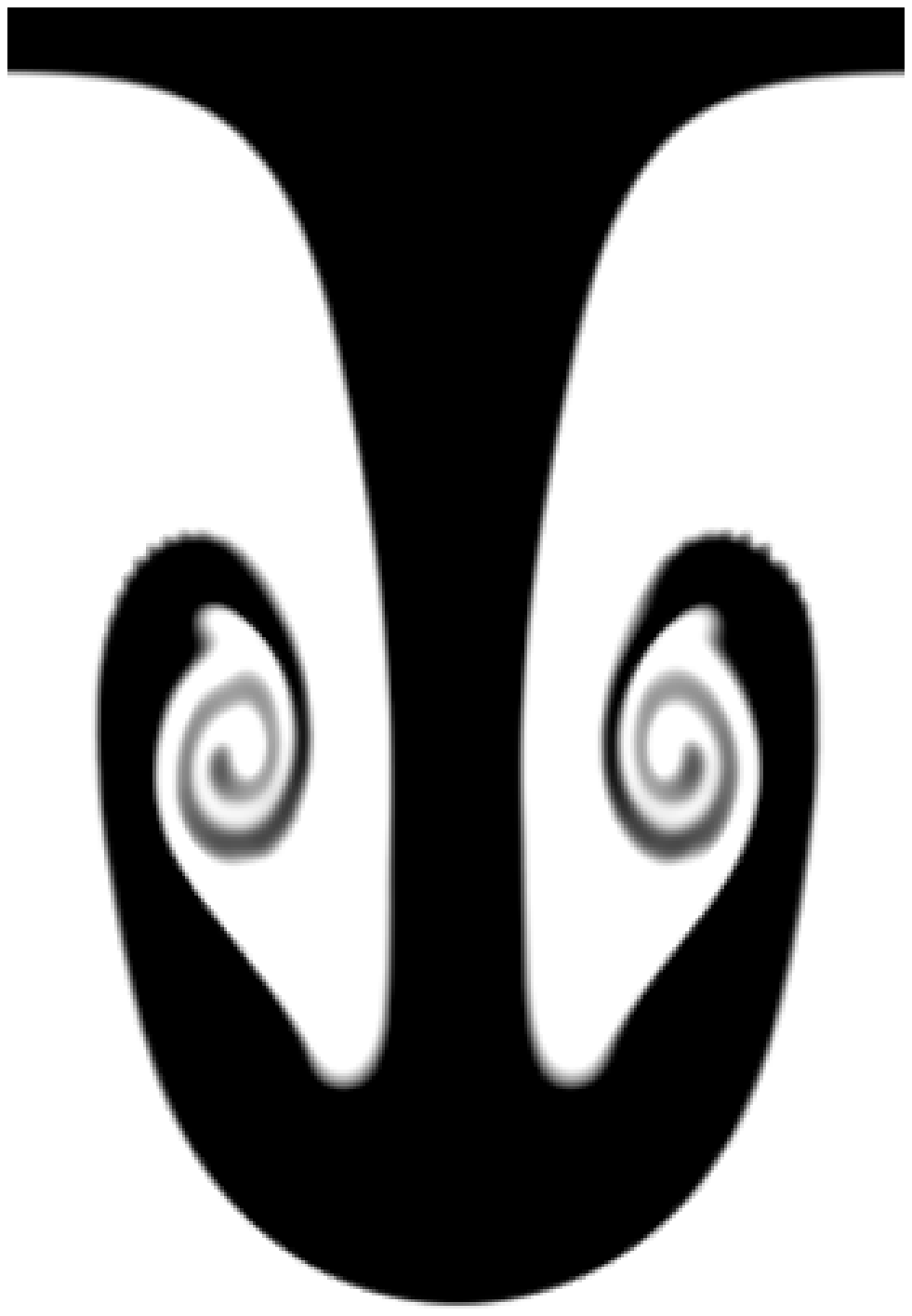}}
  \fbox{\includegraphics[width=\picscale]{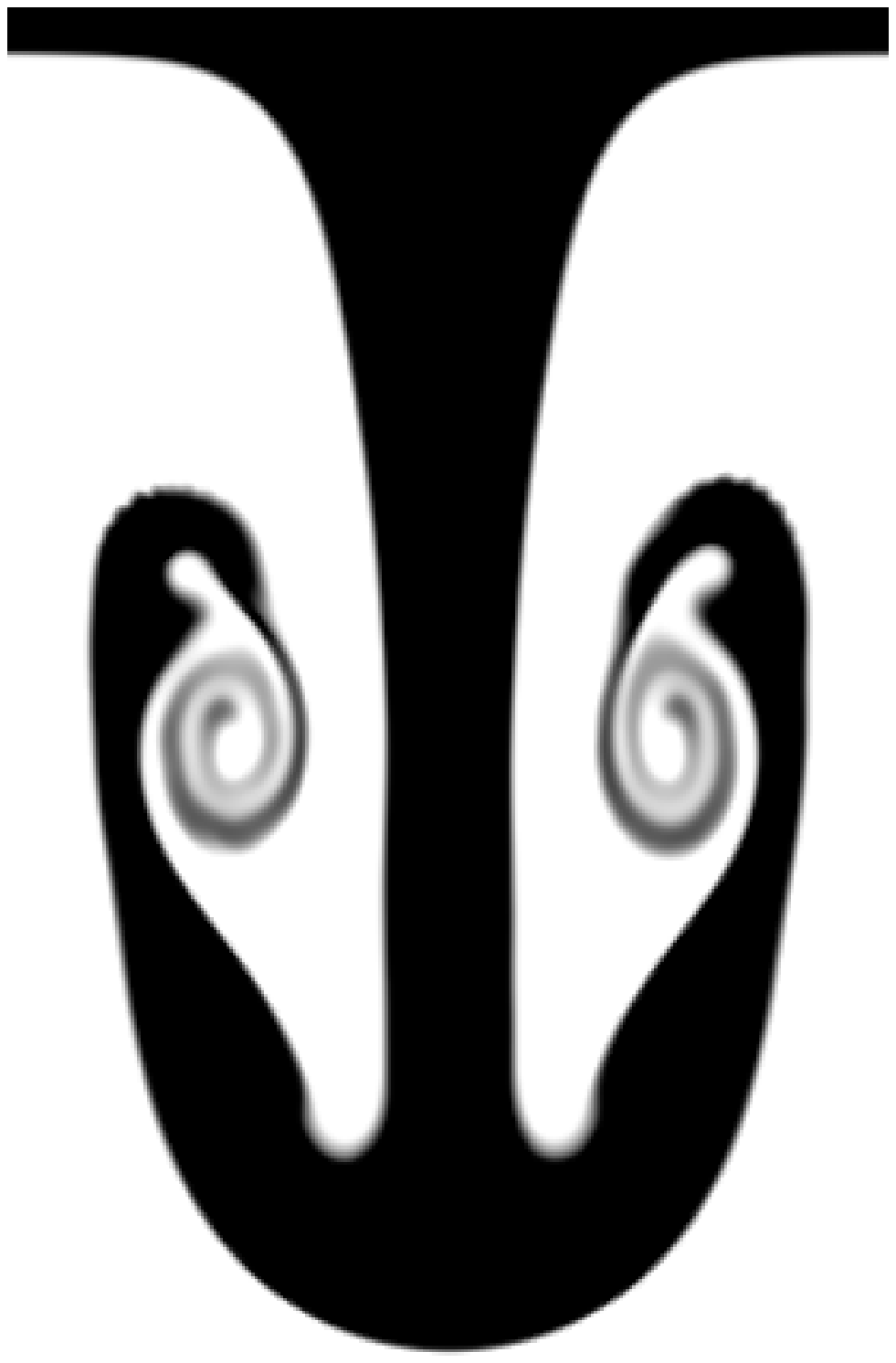}}  
  \caption{Rayleigh-Taylor instability. Atwood = 0.25 at time $T=$ 0., 0.1, 0.11, 0.12, 0.13, 0.14 (to be read linewise from top left to bottom right).}
  \label{fig:rt}
 \end{figure}

The Rayleigh-Taylor instability is probably the most famous example of
a hydrodynamic instability, first investigated by \cite{Rayleigh1883} and
\cite{Taylor1950}. A short overview about numerical approaches can be
found in \cite{Tryggvason1988}. A more recent implementation is described in
\cite{DSS07}. 
We use stabilized P1-P1 elements and test our scheme for parameters $\eta_1 = \eta_2 = 10^{-3}$, Atwood number
$A=0.25$ ($\hat\rho_\text{avg}= 0.001$), gravity $\kflow_\text{grav}=(0,-10^5)^T$, mobility $M=0.01$, $\delta=0.1$, $\sigma=0.1$. 
It is well-known that the simulation of Rayleigh-Taylor instabilities is very
demanding with respect to the discretization both in space and in time. In
Figure~\ref{fig:rt}, we find quite a good resolution of the occuring
eddies. Only at very large times -- when structures are getting finer and
finer and when mixture regions start to fatten, a loss of symmetry becomes
visible. The maximum level of adaptive spatial refinement is 16 (bisectional
grid refinement) -- this means a maximum of $1.3\cdot 10^5$ degrees of freedom.

\subsection{A first comparison with the model by \texorpdfstring{\cite{DSS07}}{Ding et al.}}\label{subsec_frame_indifferent} 
The last numerical experiment is devoted to the question to which extent the term
 $\rkla{\rkla{\rhoDerivative \jflow} \cdot \nabla} \vflow$  -- which is not
contained  in the model \cite{DSS07} -- influences the evolution.
The paper \cite{Aland2012} already compared the performance of the models
in \cite{DSS07} and in \cite{AGG} in the framework of the 
benchmark test \cite{Hysing2009} of a rising droplet in a constant gravitational field. For this
problem, the authors in \cite{Aland2012} could not find significant
improvement by the more involved model of \cite{AGG}.

Of course, both the models in \cite{DSS07} and \cite{AGG} may serve as a diffuse-interface approximation to the same sharp-interface model of two-phase flow with different mass densities. However, the question remains, whether this approximation is always of the same quality.

\begin{figure}[ht]
  \centering

  \setlength{\picscale}{0.37\textwidth}

  \tiny{AGG}
   \includegraphics[width=\picscale]{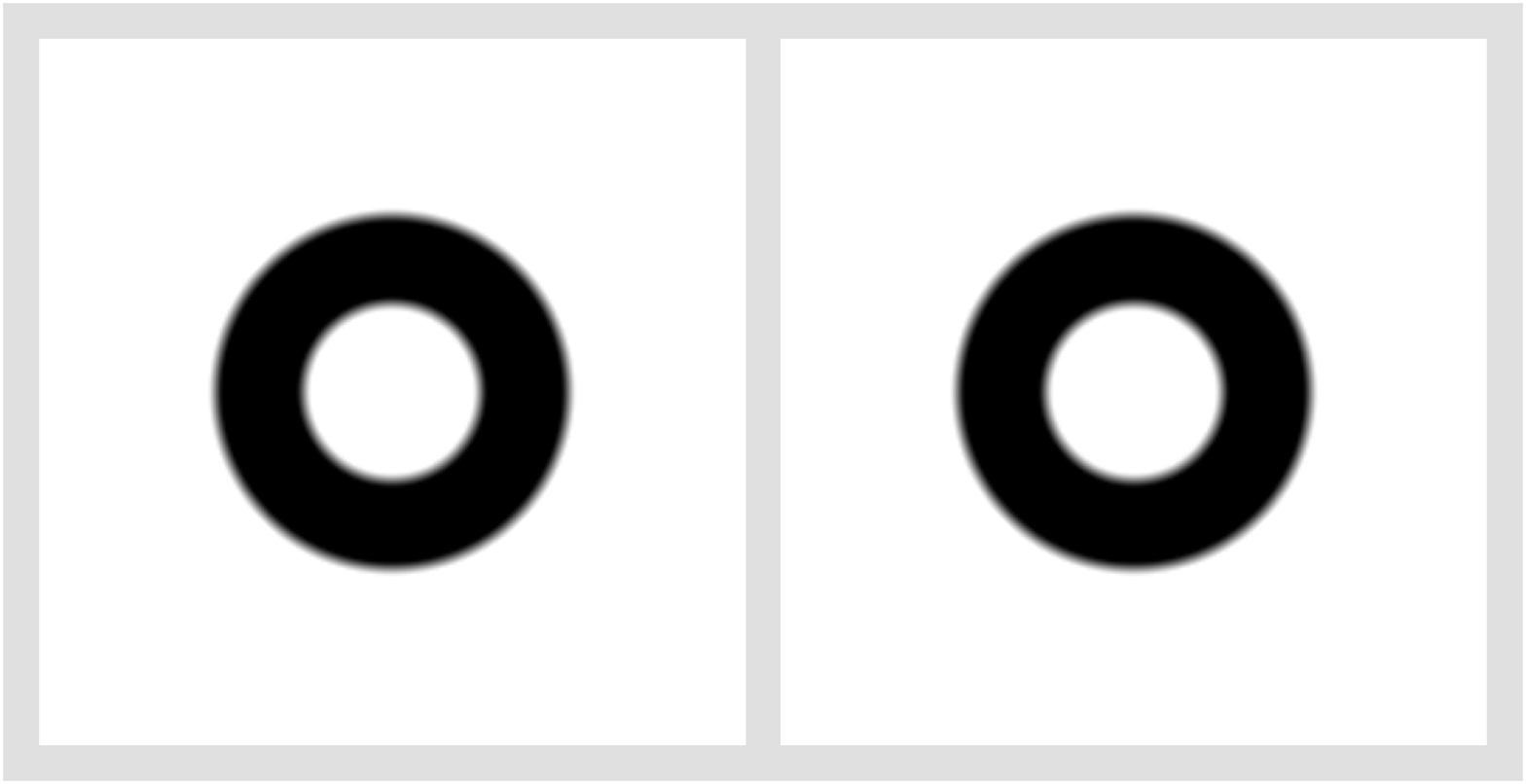}
   DSS \hfill AGG
   \includegraphics[width=\picscale]{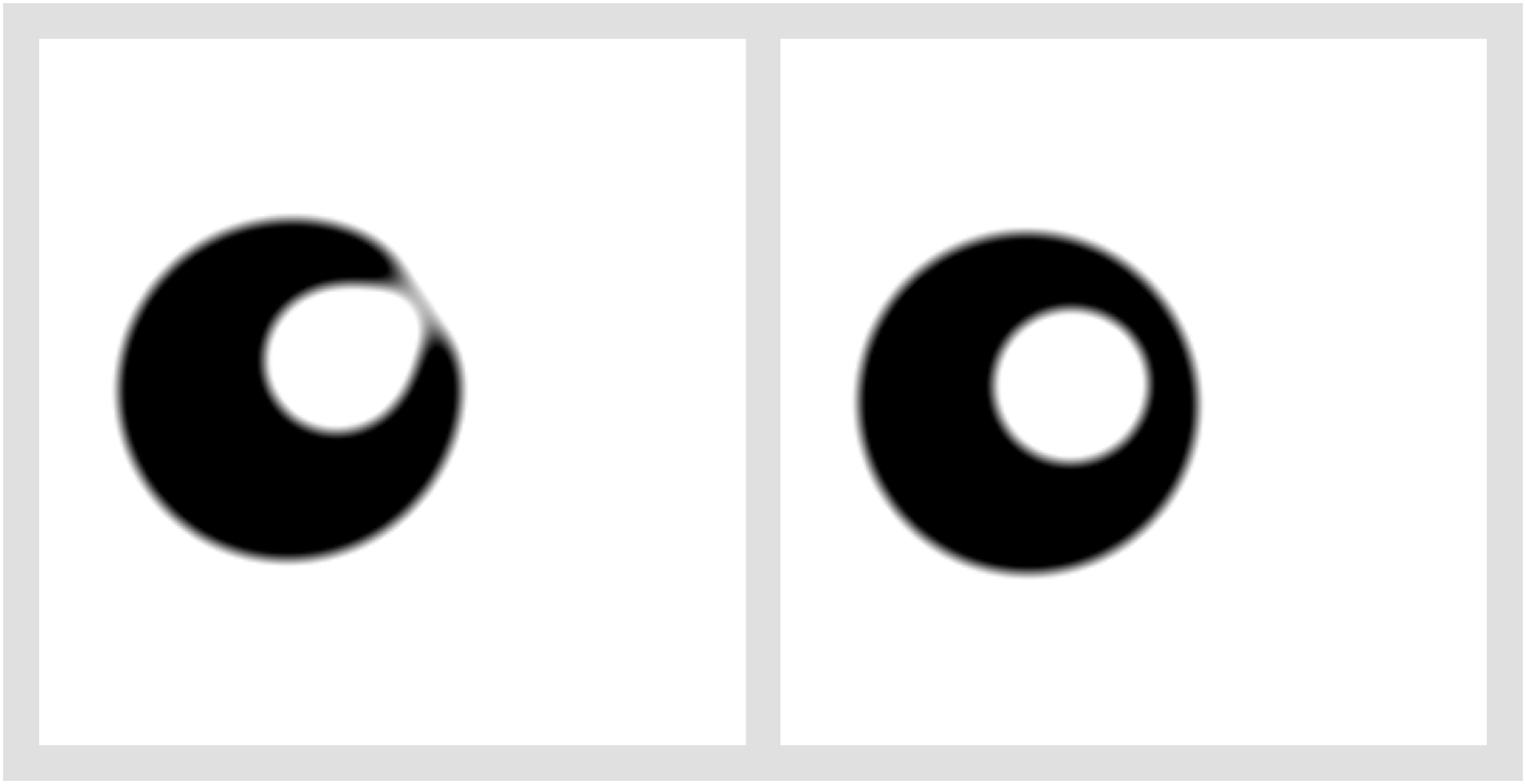}
   DSS

   AGG
   \includegraphics[width=\picscale]{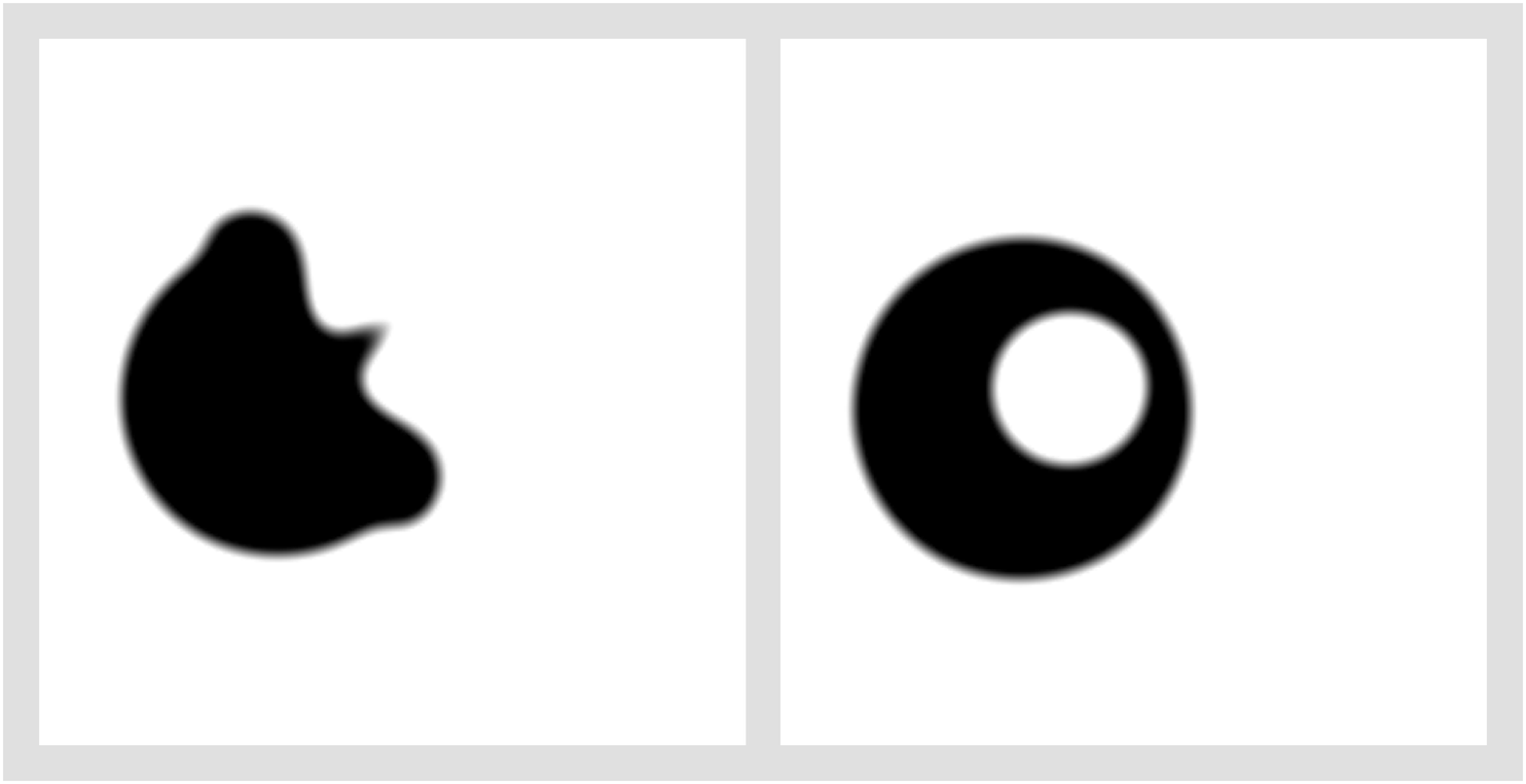}
   DSS \hfill AGG
   \includegraphics[width=\picscale]{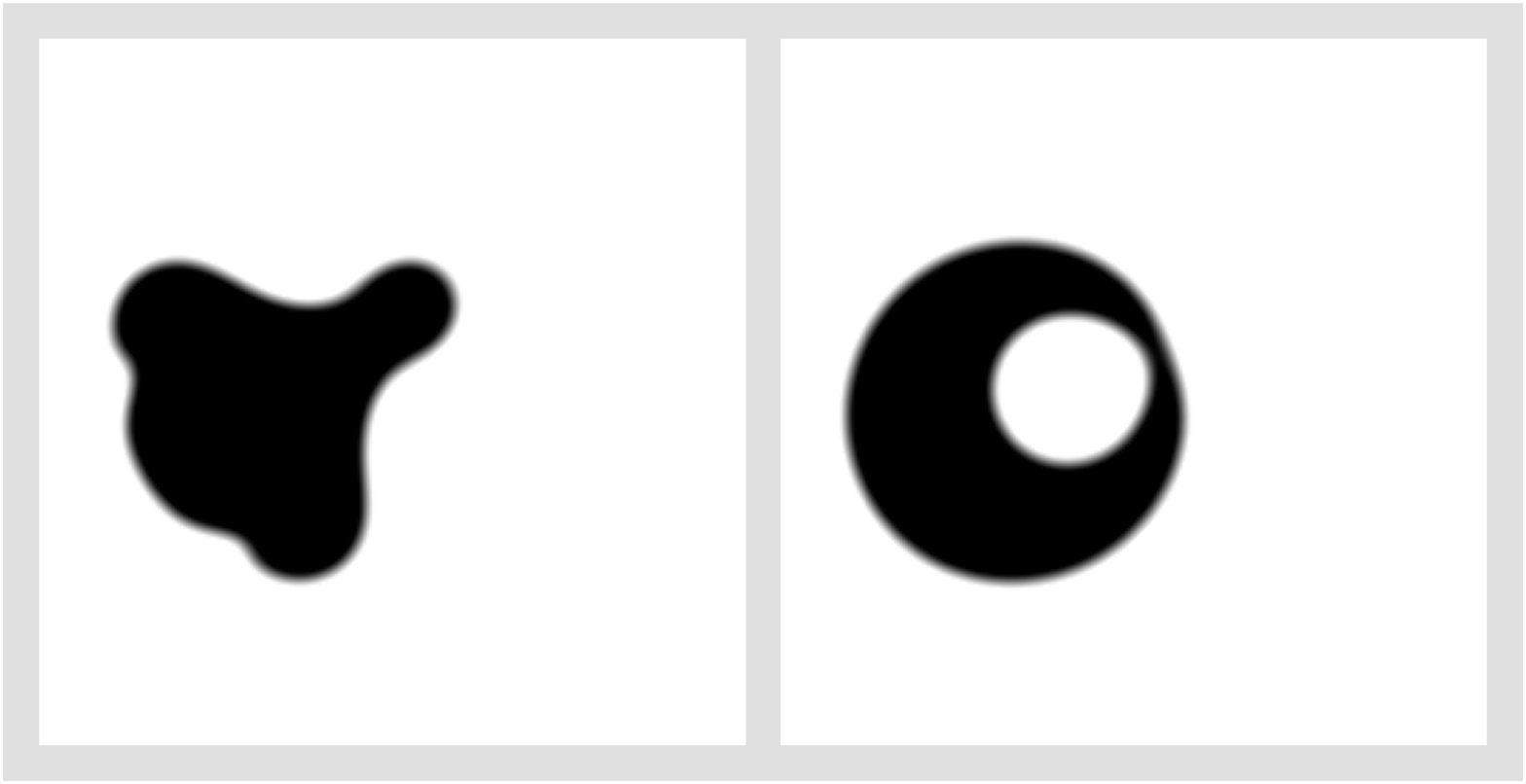}
   DSS

   AGG
   \includegraphics[width=\picscale]{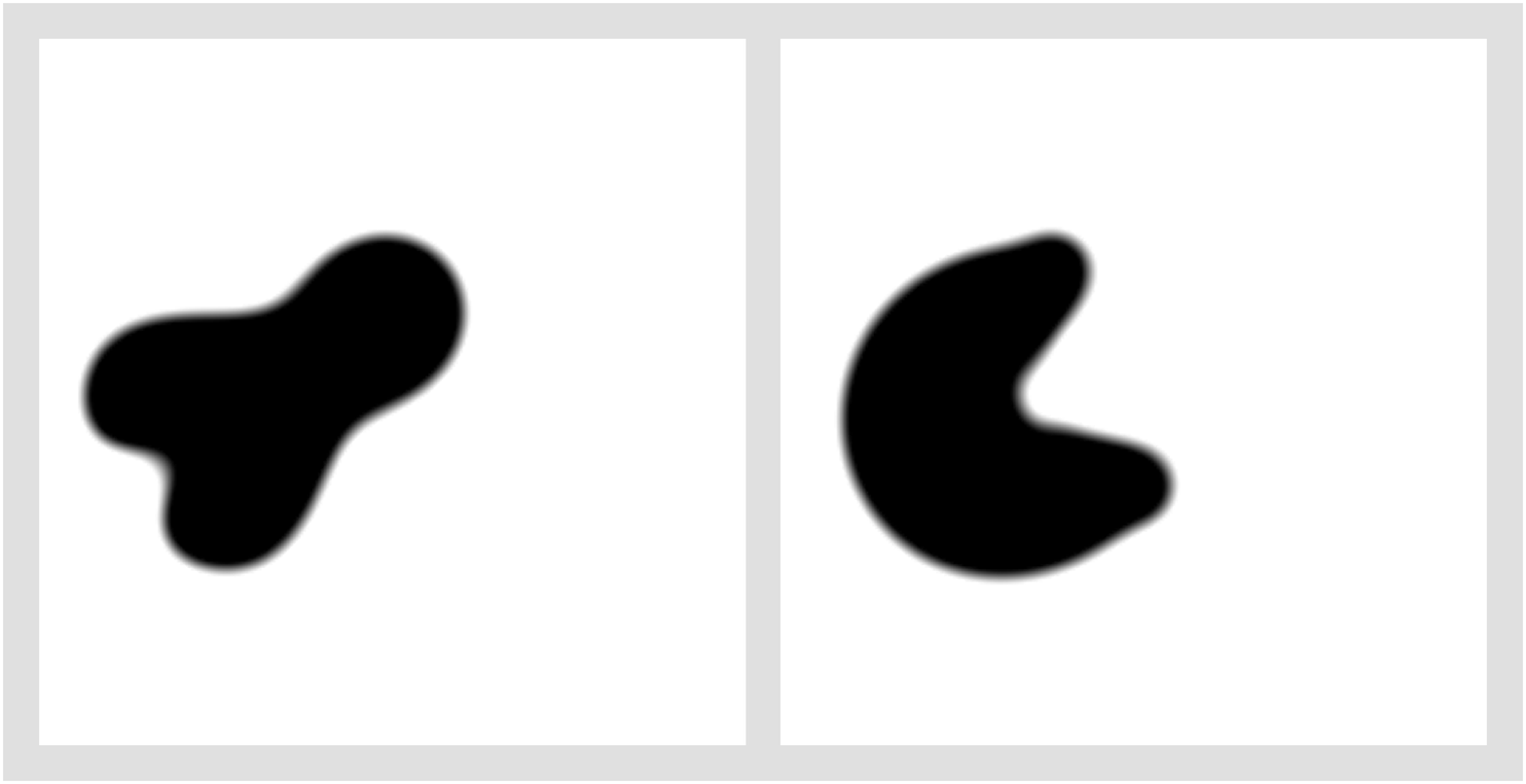}
   DSS \hfill AGG
   \includegraphics[width=\picscale]{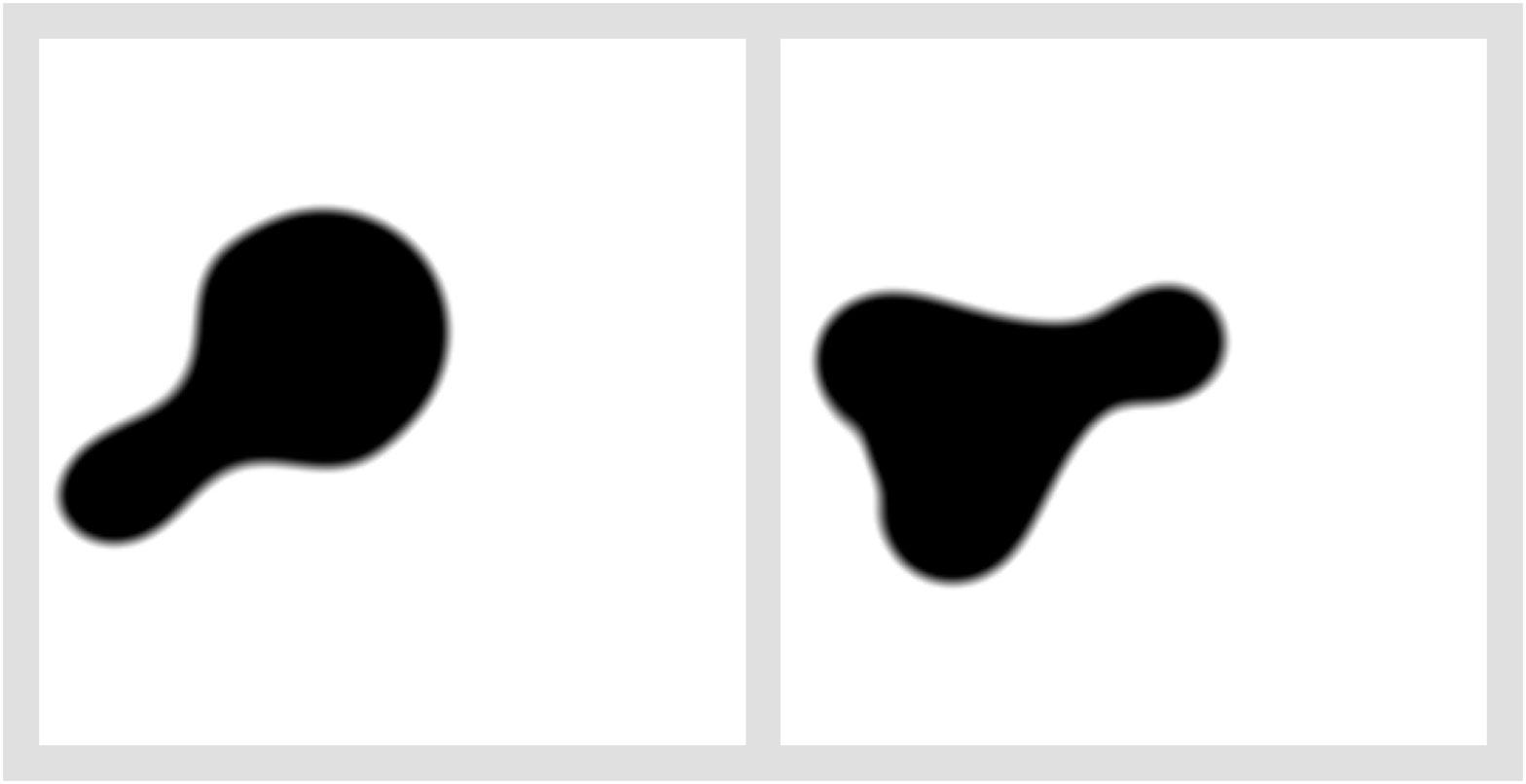}
   DSS

   AGG
   \includegraphics[width=\picscale]{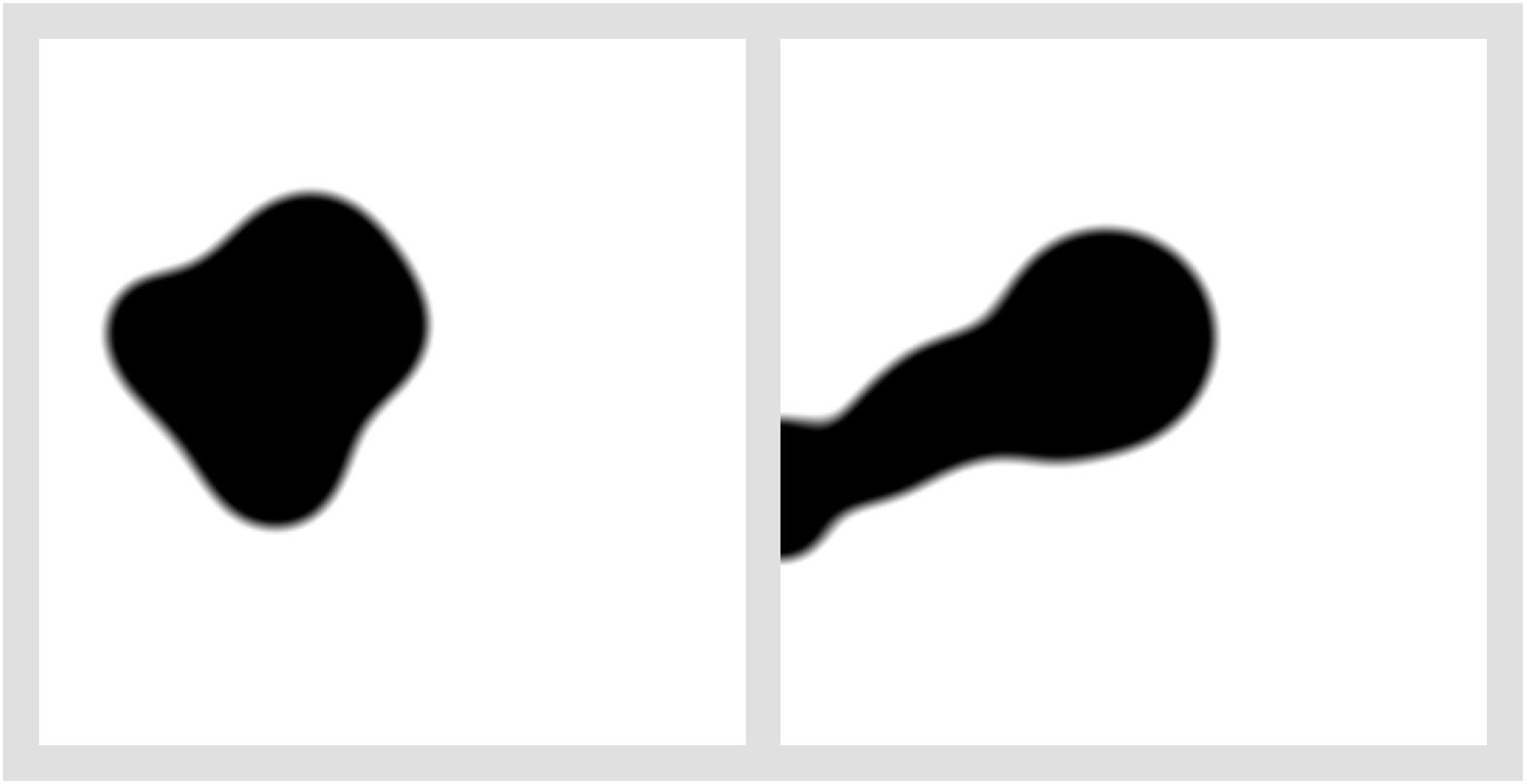}
   DSS \hfill AGG
  \includegraphics[width=\picscale]{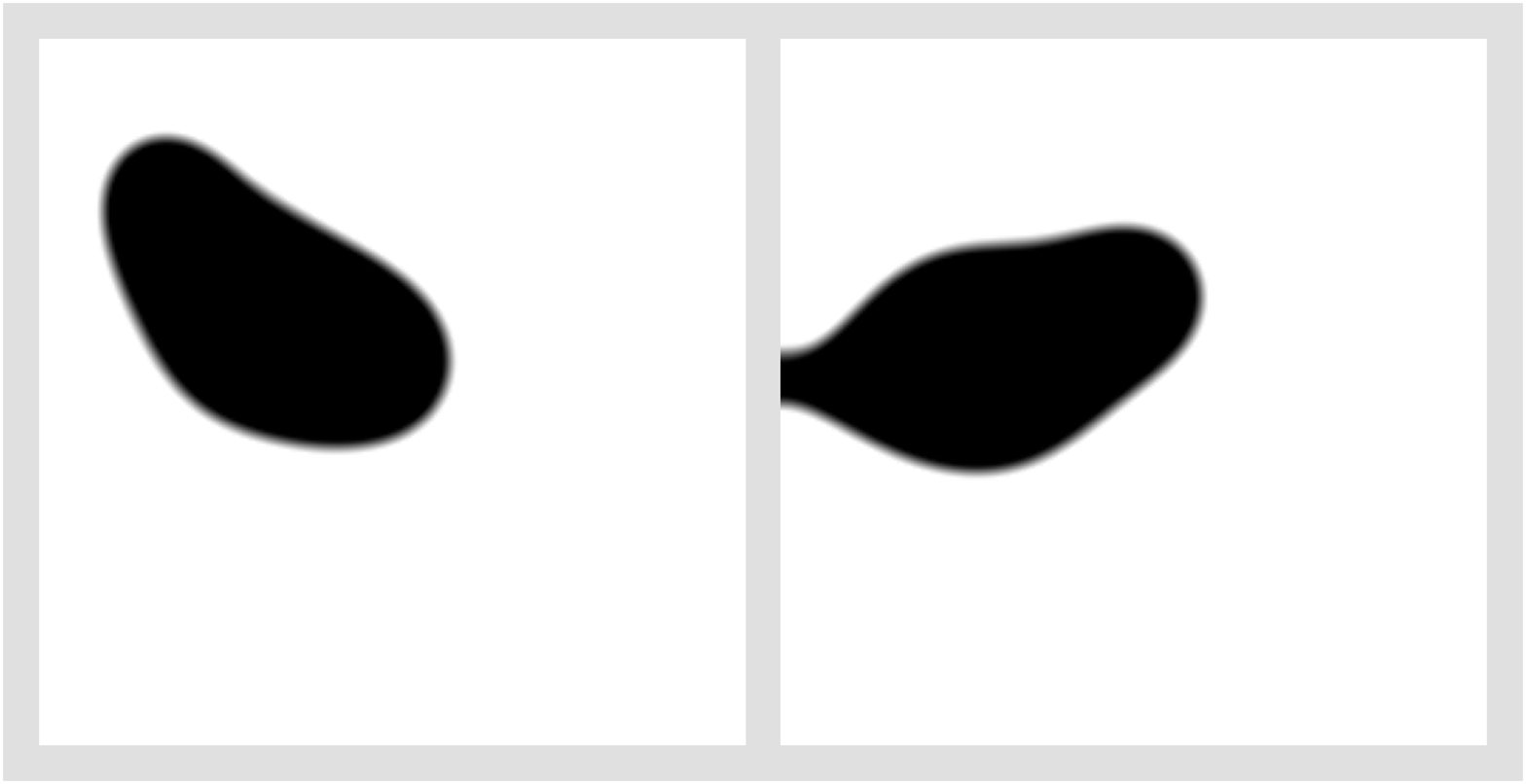}
   DSS
   
  \caption{Rotating force density field, juxtaposition of simulations based on \cite{AGG} (left) 
    and on \cite{DSS07} (right) at times $T$ = 0, 0.26, 0.27, 0.28, 0.29, 0.31, 0.33, 0.37 (to be read linewise from top left to bottom right). }
  \label{fig:frame_indifferent}
 \end{figure}

In this paper, we present simulations which show considerable differences. We set up an experiment in a rotating force field and  
we compare the results obtained by the two models. As initial datum, we take
an annular droplet ($\hat\rho_2=0.019$) surrounded by a much lighter fluid
($\hat\rho_1=0.001$). 
We assume the force density $\kflow_\text{grav}$ to rotate around the origin, i.e. 
$\kflow_\text{grav} = \mathbf O(t)\cdot (0,100)^T $ with $\mathbf O(t)$ a rotation field
with constant angular velocity of five rotations per unit of time.
Figure~\ref{fig:frame_indifferent} presents a sequence of juxtapositions of
the two experiments at selected time-steps. The right image (labeled by ``DSS'')
always corresponds to a simulation based on the model \cite{DSS07}, the left
one (labeled by ``AGG'') corresponds to a simulation based on \cite{AGG}. For the simulation of the model in \cite{DSS07}, we used the P1-P1 version of our code for \cite{AGG}, of course with the only modification, that the $\jflow$-coupling in the momentum equation was omitted.

It is worth mentioning that due to the non-constant gravitational field the symmetry is lost immediately. 
Therefore, the annulus gets perturbed and is no longer stable. Interestingly, the topological change of annulus breakup occurs earlier in the sequence based on the model \cite{AGG} than in the sequence based on \cite{DSS07}. Whether this observation can be explained by consistency with thermodynamics in the sense that energy dissipation is enhanced in \cite{AGG}, may be the subject of further studies.

\section{Conclusion}\label{sec:conclusion}
We proposed a numerical scheme for the diffuse-interface model in \cite{AGG} for two-phase flow with incompressible liquids of general mass densities. The scheme is stable in the sense that for discrete solutions discrete counterparts of the total physical energy at times $t_2>t_1$ are bounded by the sum of this energy at time $t_1$ and the work done by external forces during the time interval $(t_1, t_2)$.

We presented two different, mathematically equivalent formulations and gave details of a fully practical splitting scheme for one of them. Numerical experiments in 2D show experimental convergence of the scheme and underline its good performance, in particular in numerically demanding experiments related to tip formation, Rayleigh-Taylor instability and topology changes.

\paragraph{\bf Acknowledgement} The authors gratefully acknowledge support by the
priority programme SPP1506 ``Transport processes at fluidic interfaces'' of
German Science Foundation (DFG). In particular, Stefan Metzger, a student co-worker, contributed to the section on numerical experiments.


\bibliographystyle{amsplain}
\addcontentsline{toc}{section}{References}

\providecommand{\bysame}{\leavevmode\hbox to3em{\hrulefill}\thinspace}
\providecommand{\MR}{\relax\ifhmode\unskip\space\fi MR }
\providecommand{\MRhref}[2]{%
  \href{http://www.ams.org/mathscinet-getitem?mr=#1}{#2}
}
\providecommand{\href}[2]{#2}

\end{document}